\documentclass[11pt]{article}
\usepackage{a4,amsfonts,amssymb,amsmath}
\usepackage{epsfig} 
\usepackage{psfrag}
\usepackage{graphicx}

\def\proof{{\bf Proof:\quad}}

\def\qed{\hfill\rule{2.2mm}{2.2mm}\vspace{1ex}}

\newtheorem{theorem}{Theorem}[section]

\newtheorem{corollary}[theorem]{Corollary}

\newtheorem{lemma}[theorem]{Lemma}
\newtheorem{proposition}[theorem]{Proposition}
\newtheorem{remark}[theorem]{Remark}

\def\eps{\varepsilon}
\def\w{\omega}
\def\ohne{\backslash}
\def\CC{\mathbb C}

\def\ZZ{\mathbb Z}
\def\DD{\mathcal D}
\def\FF{\mathcal F}

\def\NN{\mathbb N}
\def\RR{\mathbb R}
\def\TT{\mathbb T}
\def\LL{\mathcal L}

\def\BB{\mathcal{B}}
\def\TT{\mathcal{T}}
\def\FF{\mathcal{F}}

\newcommand{\HH}{\mathcal H}

\newcommand{\re}{\mathop{\rm Re}\nolimits}
\renewcommand{\Re}{\re}
\newcommand{\im}{\mathop{\rm Im}\nolimits}
\renewcommand{\Im}{\im}

\def\text{\mbox}

\title{On Laplace--Carleson embedding theorems}
\author{{Birgit Jacob\thanks{Fachbereich C - Mathematik und Naturwissenschaften,
Bergische Universit\"at Wuppertal, Gau\ss stra\ss e 20, 42097 Wuppertal, Germany, \tt jacob@math.uni-wuppertal.de}
\qquad Jonathan R.~Partington\thanks{School of Mathematics,
University of Leeds,
Leeds LS2 9JT, U.K. \tt
J.R.Partington@leeds.ac.uk}
\qquad Sandra Pott\thanks{Centre for Mathematics, Faculty of Science, Lund University, S\"olvegatan 18, 22100 Lund, Sweden, \tt Sandra.Pott@math.lu.se}}}

\begin{document}
\maketitle
\begin{abstract}
This paper gives embedding theorems for a very general class of weighted Bergman spaces: the
results include a number of classical
Carleson embedding theorems as special cases.    The little Hankel operators on these Bergman spaces are also considered.                   Next, a study is made of Carleson
embeddings in the right half-plane induced by taking the Laplace transform
of functions defined on the positive half-line (these embeddings have applications in control theory):
particular attention is given to the case of a sectorial measure or a measure supported on a strip,
and complete necessary and sufficient conditions for a bounded embedding are given in many cases.
\end{abstract}


{\bf 2000 Subject Classification.} 30D55, 30E05, 47A57, 47D06, 93B05.


\section{Introduction and Notation}

This paper brings together two ideas which are fundamental to the study of Banach spaces of analytic functions. First, there is the classical Paley--Wiener theorem (see, for example \cite{nik}). This asserts that the Laplace transform (respectively, the Fourier transform) provides an isometric
isomorphism between the space $L^2(0,\infty)$ and the Hardy space $H^2$ of the right-hand half plane $\CC_+$ (respectively, the upper half plane). There are many generalizations of this result given, as mentioned below; for example the weighted space $L^2(0,\infty; dt/t)$ corresponds to the Bergman space $A^2(\CC_+)$. The strongest known generalization of this result is to the class of so-called Zen spaces (in fact we provide a slightly stronger generalization below), determined by an appropriate measure $\nu$ on the half-plane.\\

The other theme that features in this paper is the notion of a Carleson embedding. The classical Carleson theorem is to do with 
finding a simple condition for
the boundedness of the canonical injection $H^2(\CC_+) \to L^2(\mu)$, where $\mu$ is a Borel measure on $\CC_+$, and has many applications in
function theory and harmonic analysis (again, \cite{nik} is a good reference for this). The embedding theorem has also been generalized to the Bergman space.\\

Bringing these ideas together, we characterise the boundedness of Carleson embeddings for the very general class of Zen spaces (subject to
a necessary technical condition on the measure $\nu$), prove new results about  the little Hankel operators on these spaces, and discuss in addition those
embeddings induced by the Laplace transform. Ultimately, this allows us to give precise results even in the case of $L^p$ for $p \ne 2$,
when there is no exact generalization of the Paley--Wiener theorem.\\

Our results have applications in terms of interpolation in certain spaces of holomorphic functions and also admissibility and controllability in diagonal semigroups, as outlined in Section~\ref{sec:applic}. (Full details will be
presented elsewhere \cite{jpp12}).\\

So, let $w$ denote a weight function on the imaginary axis $i\RR$, let
$\mu$ be a positive regular Borel measure on the right half plane $\CC_+$ and let
$1 \le p,q \le \infty$. Embeddings of the form
\begin{equation} \label{eq:carlemb}
   L^p_w(i\RR) \rightarrow L^q(\CC_+, \mu),
\end{equation}
where a locally integrable function $f$ on the imaginary axis $i\RR$
is mapped to its Poisson extension on the right half plane $\CC_+$,
are known as Carleson  embeddings, and have been much
studied in the literature. In linear control, another, related
class of embeddings plays an important role, namely embeddings of the
form
$$
  \HH^p_{\beta, w}(0,\infty) \rightarrow L^q(\CC_+, \mu), \quad f
  \mapsto \LL f= \int_0^\infty e^{-t \cdot} f(t) dt,
$$
given by the Laplace transform $\LL$. Here,  $\HH^p_{\beta, w}$
denotes the Sobolev space of index $\beta$ and weight $w$: the case $\beta=0$ corresponds to
a weighted $L^p$ space. We shall
refer to such embeddings as Laplace--Carleson embeddings.



In the easiest case, $\beta=0$, $w \equiv 1$ and $p=q=2$, the Laplace
transform maps
$\HH^2_{0,1}(\RR_+)=L^2(\RR_+)$ isometrically up to a constant to
$H^2(\CC_+)$, which is a closed subspace of $L^2(i\RR)$, and we 
only have to deal with the unweighted classical Carleson embedding theorem
for $p=2$. This can be found in many places, for example, \cite{garnett,nik}. 

\begin{theorem}[Carleson embedding theorem] \label{thm:ccarleson}
Let $\mu$ be a positive regular Borel measure on the  right half plane
$\CC_+$. Then the following are equivalent:
\begin{enumerate}
\item The natural embedding
$$
     H^p(\CC_+) \rightarrow L^p( \CC_+, \mu)
$$
is bounded for some (or equivalently, for all) $1 \le p < \infty$.
\item
There exists a constant $C> 0$ such that
$$
    \int_{\CC_+} | k_\lambda(z)|^2 d \mu(z) \le C \|k_\lambda\|_{H^2}^2 \text{ for all } \lambda \in \CC_+,
$$
where $k_\lambda(z)=\frac{1}{2\pi} \frac{1}{z+\overline\lambda}$
for $\lambda, z \in \CC_+$.
\item
$$
    \mu(Q_I) \le C |I| \text{ for all intervals } I \subset i\RR,
$$
where $Q_I$ denotes  the {\em Carleson square} $Q_I = \{ z= x + iy \in \CC_+: iy \in I, 0 < x < |I|\}$.
\end{enumerate}
In this case, $\mu$ is called a {\em Carleson measure}.
\end{theorem}

A further relatively easy case is $p=q=2$ and $w$ a power weight, $w(t) = t^\alpha$ with $\alpha <0$. This case
corresponds to the classical embedding theorem for standard weighted Bergman spaces on the half plane by Duren, see
e.g. \cite{duren}. For more general weights $w$, the Laplace--Carleson embedding corresponds to a new embedding theorem for weighted Bergman spaces $A^2_\nu$
on the half plane with a translation-invariant measure $\nu$, which is the subject of Section \ref{sec:zen}.

In the case of general $1\le p,q \le \infty$, $p>2$, the Laplace--Carleson embeddings are very subtle even in case that $w=1$,
due to the oscillatory part of the Laplace transform integral
kernel. A general characterization seems out of reach at the moment, but with additional conditions on the support of the measure
$\mu$, a full characterization can sometimes be given. This is the content of Section \ref{sec:pq}.

Some applications are outlined in Section~\ref{sec:applic}.\\



The reproducing kernel functions for $H^2(\CC_+)$ are denoted by $k_\lambda$, $\lambda \in \CC_+$, where $k_\lambda(z)=\frac{1}{2\pi} \frac{1}{z+\overline\lambda}$
for $z \in \CC_+$, and satisfy $f(\lambda)=\langle f, k_\lambda \rangle$ for $f \in H^2(\CC_+)$. Note that $\|k_\lambda\|^2=\frac{1}{4\pi\re\lambda}$. We will also frequently use the Poisson kernel
$$
    p_\lambda(t)=\frac{1}{\pi} \frac{y}{x^2+(y-t)^2}, \qquad (z=x+iy \in \CC_+, \  t \in \RR),
$$
which is related to the reproducing kernel $k_\lambda$ by
$$
   \frac{1}{\|k_\lambda\|^2} |k_\lambda(i t)|^2 =    p_\lambda(t)   \quad \text{ for }  \quad   t \in \RR.
$$

\section{Embedding theorems for weighted Bergman spaces}   \label{sec:zen}
In this section, we will be interested in embeddings
$$
   A_\nu^p(\CC_+) \hookrightarrow L^p(\CC_+, \mu),
$$
where  $A_\nu^p(\CC_+)$ is a  weighted Bergman space defined below and
$\nu$ is a translation-invariant positive regular Borel measure on $\overline{\CC_+}$, that is, 
$\nu =  \tilde \nu \otimes  \lambda$, where $\lambda$ denotes Lebesgue measure and $\tilde \nu$ is a positive  regular Borel measure on $[0, \infty)$. This
corresponds to the case of radial measures on the unit disc.

The investigation of such embeddings has a long history, starting with \cite{duren1}  for the case of the standard weighted Bergman space on $\CC_+$ with
$d \tilde \nu (t) = t^\beta dt$ for $\beta > -1$, respectively the standard weights $(1 -|z|)^\beta$ on the disc. 
Oleinik  \cite{ol} observed already in 1974
that for measures $\rho(t)dt$, where the weight $\rho$ decreases very fast towards $t=0$, such as $\rho(t) = e^{-1/t^{1 + \gamma}}$, $\gamma>0$, 
it is not sufficient to compare the weights of Carleson squares
$\mu(Q_I)$ and $\nu(Q_I)$ (or equivalently, to compare the measures of Euclidean balls $D$ centered on the imaginary axis). Instead, in the example above
one has to consider the measures of 
Euclidean balls away from the imaginary axis, 
$$
     D_z = D \left(z,   \frac{(\Re z)^{1+\gamma} }{ (1 + \Re z)^{\gamma}} \right)
$$
for $z \in \CC_+$.
Roughly speaking, the faster the weight $\rho(t)$ decreases for $t \to 0$, the smaller the radius of the ball $D_z$ in relation to the distance of $z$ to the imaginary axis, and the more detailed information
on the measure $\mu$ is required. Recently, necessary and sufficient conditions have been found for the case even faster decreasing weights $\rho(t)$, such as double exponentials
\cite{pau}. Our aim in this section is somewhat different: we want to find  a class of measures $\tilde \nu$ as large as possible, for which a characterisation
in terms of Carleson squares, and in terms of testing on powers of reproducing kernels $k_z$, $z\in \CC_+$, is possible. 

Our approach relies on a dyadic decomposition of the half-plane adapted to the measure which is somewhat technical, but has the advantage that no smoothness or continuity properties of the measure are
required. Clearly, we need a growth condition on $\tilde \nu$ in $0$ and the most natural condition to impose is the well-known $(\Delta_2)$
(doubling) condition at $0$. 

One motivation  to treat this general setting, apart from the interest of Carleson-type embedding theorems in their own right, is to obtain Laplace-Carleson embedding theorems on a large class 
class of weighted spaces $L^2_w(0, \infty) $  (see Theorem \ref{thm:CarlLap} below), which are important in Control Theory. This is exploited in \cite{jpp12}. 
On the other hand, the general embedding theorem makes it also possible to study other
important operators of analytic function spaces, such as Hankel operators and Volterra-type integration operators, in a more general setting. Such operators on weighted Bergman spaces for various classes
of weights have for
example been studied in \cite{rochberg}, \cite{aleman} and recently in \cite{pau}, \cite{pau1}. As an example of this type of application, we study little Hankel operators on Zen spaces in Subsection \ref{sec:2.2}.

\subsection{Carleson measure on Zen spaces}

Let $\tilde \nu$ be a positive  regular Borel measure on $[0, \infty)$ satisfying the following $(\Delta_2)$-condition:
\begin{equation}\tag{$\Delta_2$}
   R:=  \sup_{t >0} \frac{\tilde \nu[0, 2t)}{\tilde \nu[0, t)} < \infty.
\end{equation}

This is sometimes referred to as a {\em doubling condition}, and such measures have been studied in the theory of harmonic analysis and partial differential equations for many years (an early reference is \cite{ST80}).
Let $\nu$ be the positive regular Borel measure on $\CC_+ =[0, \infty) \times \RR $ given by $d\nu = d \tilde \nu \otimes d \lambda$, where $\lambda$ denotes Lebesgue measure.
In this case, for $1\le p<\infty$, we call 
$$A^p_\nu= \left\{ f: \CC_+ \rightarrow \CC \text{ analytic}: \sup_{\eps >0} \int_{\overline{\CC_+}} |f(z+\eps)|^p d\nu(z)< \infty \right\}
$$ a {\em Zen space on $\CC_+$}. If $\tilde \nu(\{0\}) >0$, then by standard Hardy space theory, $f $ has a well-defined boundary function $\tilde f \in L^p(i \RR)$, and we can give meaning
to the expression $\int_{\overline{\CC_+}} |f(z)|^p d\nu(z)$.   Therefore, we may write
$$
                         \|f\|_{A^p_\nu} = \left(\int_{\overline{\CC_+}} |f(z)|^p d\nu(z) \right)^{1/p}.
$$
Note that this expression makes sense in the case that $\tilde \nu(\{0\}) =0$ (e.g. the Bergman space), since $f$ is still defined $\nu$-a.e.\ on $\overline{\CC_+}$.
Clearly the space $A^2_\nu$ is a Hilbert space.

Well-known examples of Zen spaces are  Hardy space $H^p(\CC_+)$, where $\tilde \nu$ is the Dirac measure in $0$, or the standard weighted Bergman spaces $A^p_{\alpha}$, where 
$d \tilde \nu(t) = t^\alpha dt $, $ \alpha > -1$.  
Some further examples constructed from Hardy spaces on shifted half planes
were given by Zen Harper in \cite{zen09,zen10}.
Note that by the $(\Delta_2)$-condition, there exists $N \in \NN$ such that $k_\lambda^N \in A^p_\nu$ for all $\lambda \in \CC_+$ and all $1 \le p < \infty$.
Here is our Embedding Theorem for Zen spaces.
\begin{theorem}   \label{thm:Zen}
Let $1 \le p < \infty$, let $A^p_\nu$ be a Zen space on $\CC_+$, with measure $\nu = \tilde \nu \otimes \lambda$ as above, and let $\mu$ be a positive regular Borel measure on $\CC_+$. 
Then the following are equivalent:
\begin{enumerate}
\item The embedding $A^p_\nu \hookrightarrow  L^p(\mathbb C_+,\mu)$ is well-defined and bounded for one, or equivalently for all, $1 \le p < \infty$.
\item For one, or equivalently for all, $1 \le p < \infty$, and some sufficiently large $N \in \NN$, there exists a constant $C_p>0$ such that
\begin{equation}   \label{eq:zenreprocond}
    \int_{\CC_+}    |(k_{\lambda}(z))^N|^p d \mu(z)  \le C_p     \int_{\CC_+}    |(k_{\lambda}(z))^N|^p  d \nu(z)       \text{ for each } \lambda \in \CC_+.
 \end{equation}

 \item There exists a constant $C>0$ such that 
\begin{equation}   \label{eq:zencond}
 \mu(Q_I) \le C \nu(Q_I) \text{ for each Carleson square }Q_I.
 \end{equation}
\end{enumerate}
\end{theorem}
\proof
The implication (1) $\Rightarrow$ (2) for fixed $p$ is immediate. 

For the implication (2) $\Rightarrow$ (3), we use a standard argument using the decay of reproducing kernels. 
Given an interval $I$ in $i \RR$, and let $\lambda$ denote the centre of the Carleson square $Q_I$.
Note that
$$
    |(k_\lambda(z))^N|  \ge \frac{1}{(2\pi)^N (4 \Re \lambda)^N}   \text{ for } z \in Q_I.
$$
Hence 
$$
      \int_{\CC_+}    |(k_{\lambda}(z))^N|^p d \mu(z)   \ge \frac{1}{(8\pi)^{pN} ( \Re \lambda)^{pN}}  \mu(Q_I) .
$$
It only remains to estimate  $\int_{\CC_+}    |(k_{\lambda}(z))^N|^p  d \nu(z)$ in terms of $\nu(Q_I)$. 
Let $R$ be the constant from the $(\Delta_2)$-condition. 

For $k \in \NN$, let $2^k I$ denote the interval with the same centre as $I$ and the $2^k$ fold length.  By  the $(\Delta_2)$-condition, 
$\nu(Q_{2^k I}) \le R^{k} 2^k \nu(Q_I)$. Note that
$$
     |(k_\lambda(z))^N|  \le      \frac{1}{(2\pi)^N}     \frac{1}{ (2^{k-1}  \Re \lambda)^N  }   \text{ for  } z \in  Q_{2^k I} \backslash  Q_{2^{k-1} I}.
$$

 Hence
\begin{multline*} 
               \|k_{\lambda}^N\|^p_{A_\nu^p}  \\
                  \le \sum_{k=0}^\infty           \nu(Q_{2^k I})     \left(\frac{1}{(2\pi)^N}     \frac{1}{ (2^{k-1}  \Re \lambda)^N  } \right)^p               
                                  \le\nu(Q_I)           \frac{1}{(2\pi)^{pN}}   \frac{1}{    (\Re \lambda)^{Np} }     \sum_{k=0}^\infty    \frac{2^k R^{k}   }{2^{(k-1)Np}}.  
\end{multline*}
Choosing $N$ sufficiently large, depending on $R$,
we find that the sum on the right converges to a constant $K_{N,p}$ depending on $N$ and $p$. Hence
\begin{multline}
         \frac{1}{(8\pi)^{pN} ( \Re \lambda)^{pN}}      \mu(Q_I)    \le  \int_{\CC_+}    |(k_{\lambda}(z))^N|^p d \mu(z) \le C_p  \int_{\CC_+}    |(k_{\lambda}(z))^N|^p d \nu(z)  \\
            =     C_p K_{N,p}      \frac{1}{(2\pi)^{pN}}   \frac{1}{    (\Re \lambda)^{Np} }  \nu(Q_I) , 
 \end{multline}
and we obtain $ \mu(Q_I) \le C_{N,p}  \nu(Q_I)$, with a constant $C_{N,p}$  depending only on $N$ and $p$ (and hence on the $(\Delta_2)$-condition constant $R$).

(3) $\Rightarrow$ (1):

Our strategy is to deduce the boundedness of the embedding from the classical Carleson Embedding Theorem via a suitable decomposition.

Suppose for the moment that $\tilde \nu (\{0\}) =0$ and that there exists a   strictly  increasing sequence   $(a_n)_{n \in \ZZ}$  in $\RR_+$ such that
\begin{enumerate}
\item   There exists $1> c>0$ with
\begin{equation}   \label{eq:delta2}
            \frac{a_{n+1} - a_n }{ a_{n+1}} \ge c   \text{ for all } n \in \ZZ ;
\end{equation}
\item
\begin{equation}   \label{eq:delta2s}
   (2R)^3   \tilde \nu([a_{n-1}, a_{n}))   \ge     \tilde \nu([a_n, a_{n+1}))    \ge 2R  \tilde \nu([a_{n-1}, a_{n})) .
\end{equation}
\end{enumerate}
Write $\beta_n = \tilde \nu([a_n, a_{n+1})) $. Notice that by (\ref{eq:delta2s}),
\begin{multline*}
        \sum_{n \in \ZZ}    \beta_n   \|f\|^p_{H^p_{a_{n+1}}}    \le  \int_{\CC_+} |f(z)|^p d\nu   = \sum_{n \in \ZZ}  \int_{-\infty}^\infty  \int_{[a_n, a_{n+1})}   |f(t + i \w)|^p  d \tilde \nu(t) d \w  \\
             \le  \sum_{n \in \ZZ}    \beta_n   \|f\|^p_{H^p_{a_n}}      \le (2R)^3  \sum_{n \in \ZZ}    \beta_n   \|f\|^p_{H^p_{a_{n+1}}}.
\end{multline*}
Thus both  $(\sum_{n \in \ZZ}    \beta_n   \|f\|^p_{H^p_{a_{n+1}}})^{1/p} $ and  $(\sum_{n \in \ZZ}    \beta_n   \|f\|^p_{H^p_{a_n}} )^{1/p}$ give us equivalent expressions for the $A^p_\nu$-norm.
In addition, notice that if $\mu$, $\nu$ satisfy the Carleson-type Condition (\ref{eq:zencond}), then also $\mu$, $\sum_{n \in \ZZ}   \beta_n \delta_{a_n} \otimes \lambda$, satisfy the condition with the same constant $C$. Therefore, we assume without loss of generality that
$$
            \tilde  \nu = \sum_{n \in \ZZ}   \beta_n \delta_{a_n},
$$
where $\delta_{a_n}$ denotes the Dirac measure at $a_n$. For simplicity of notation we will also assume that the constant $C$ in the condition (\ref{eq:zencond}) equals $1$.

Our next step is a decomposition of $\mu$ into $\sum_{n \in \ZZ} \mu_n$, where each $\mu_n$ satisfies the Carleson-type condition  (\ref{eq:zencond}) with respect to $\nu_n = \beta_n \delta_{a_n} \otimes \lambda$.
\begin{lemma}       \label{lemm:mainest}
Let $N \in \ZZ$ and suppose that $\mu$ is supported on the closed half-plane $\overline{\CC_{a_N}}$. Then there exist positive regular Borel measures $\mu_n$, $n \ge N$, such that
\begin{enumerate}
\item
\begin{equation}    \label{cond:sum}
   \mu = \sum_{n=N}^\infty  \mu_n;
\end{equation}
\item   
\begin{equation}    \label{cond:supp}
\mu_n\text{ is supported on the closed half-plane } \overline{\CC_{a_n}};
\end{equation}
\item    There exists a constant $C' >0$ such that for all intervals $I \subset i \RR$,
\begin{equation}    \label{cond:ncarl}
   \mu_n(Q_I) \le C' \nu_n(Q_I) \quad (n  > N),  \quad   \mu_N(Q_I) \le C'  \sum_{k= -\infty}^N  \nu_k(Q_I).
\end{equation}
\end{enumerate}

Moreover, $\mu_n$ is a Carleson measure for the shifted half plane $\CC_{a_{n-1}}$, with Carleson constant
$$
         \mathcal{C}_{a_{n-1}} (\mu_n)   \le \frac{C'}{c}  \beta_n  \quad (n  > N) , \quad  \mathcal{C}_{a_{N-1}} (\mu_N)   \le \frac{C'}{c}  \sum_{k= -\infty}^N \beta_k,
$$
where $c$ is the constant appearing in the definition of the sequence $(a_n)$ above.
\end{lemma}
{\bf Proof of the lemma:}\quad  
First, we prove that $\mu_n$, $n \ge N$, exist, satisfying Conditions (\ref{cond:sum}), (\ref{cond:supp}) and (\ref{cond:ncarl}).
 By replacing $\nu_N$ with $(\sum_{n= -\infty}^N  \beta_n )\delta_{a_N}  \otimes  \lambda $,
we can assume without loss of generality that $ \nu$ is supported on $\overline{\CC_{a_N}}$.

We begin by constructing  a family of Carleson rectangles, on which (\ref{cond:ncarl}) can be checked.  Fix a dyadic grid $\DD_N$ of half-open intervals in $i \RR$ with
minimal intervals of length $a_{N}$, and denote these   intervals on length $a_{N}$ as intervals of generation $0$. The remaining intervals in the  family $\DD_N$ will be parents, grandparents, etc of
the minimal intervals.

We can assume without loss that 
\begin{equation}  \label{eq:c}
     c \le 1- \frac{1}{\sqrt{2}}.
\end{equation}
 By (\ref{eq:delta2}),
$$
          a_{N+k+1}   \ge \frac{ a_{N+k}}{1-c}  
$$
for any $k \ge 0$. If
$$
    a_{N+k+1}   \ge \frac{ a_{N+k}}{(1-c)^2} ,
$$
we choose an integer $l \ge 2$ such that
$$
     \frac{1}{1-c} \le    \gamma =       \left ( \frac{a_{N+k+1}}{a_{N+k}} \right )^{1/l}   <  \frac{1}{(1-c)^2}
$$
and add $l-1$ intermediate points 
$$
        \gamma a_{N+k}, \dots , \gamma^{l-1} a_{N+k}
$$
between $a_{N+k}$ and $a_{N+k+1}$. In this way, we create a strictly increasing sequence $(b_j)_{j \ge 0}$   such that  $b_0= a_N$, all terms of the sequence $(a_n)_{n \ge N}$ appear as terms 
of the sequence $(b_j)_{j \ge 0}$, and
\begin{equation}   \label{eq:geom}
            \frac{1}{1-c}   b_{j}   \le            b_{j+1}     <  \frac{1}{(1-c)^2} b_{j}  \text{ for all } j.
\end{equation}
The family $\FF$ of Carleson squares we want to consider is formed by rectangles of the form $ (0, b_{j+1}) \times I$, where $I \in \DD_N$, $j \ge 0$,
 for which the eccentricity is bounded above and below by
\begin{equation}   \label{eq:ecc}
                                 \frac{1}{\sqrt{2}}  <               \frac{b_{j+1}}{|I|}   \le  \sqrt{2} .
\end{equation}
Such a rectangle is denoted by $Q_{I, j}$. Note that by (\ref{eq:geom}) and (\ref{eq:c}),
 for each $I \in \DD_N$ there exits $j \ge 0$ with (\ref{eq:ecc}), and for each $j \ge 0$ there
exists exactly one size of intervals $I$ in $\DD_N$ such that  (\ref{eq:ecc}) holds.
Note that different Carleson squares in this family $\FF$  can have the same base $I  \in \DD_N$. It is easy to see that any Carleson square $Q$ over an interval in 
$i \RR$ can be covered by a bounded number of elements in $\FF$, with comparable base length. 
Therefore, up to a possible change of constant, it is sufficient to check the Carleson-type 
conditions (\ref{eq:zencond}) and (\ref{cond:ncarl}) 
on the elements of $\FF$.

\begin{figure}[t]
    \centering
   \psfrag{$b_1$}[b][b]{$b_1$}
   \psfrag{$b_2$}[b][b]{$b_2$}
   \psfrag{$b_3$}[b][b]{$b_3$}
   \psfrag{$b_4$}[b][b]{$b_4$}
   \psfrag{$b_5$}[b][b]{$b_5$}
    \psfrag{$a_N=b_0$}[b][b]{$a_N\!=\!b_0$}
   \psfrag{$a_N$}[b][b]{$a_N$}
   \psfrag{0}[t][t]{$T_{\cdot,0}$}
  \psfrag{1}[t][t]{$T_{\cdot,1}$}
  \psfrag{2}[t][t]{$T_{\cdot,2}$}
  \psfrag{3}[t][t]{$T_{\cdot,3}$}
  \psfrag{4}[t][t]{$T_{\cdot,4}$}
   \includegraphics[scale=0.55]{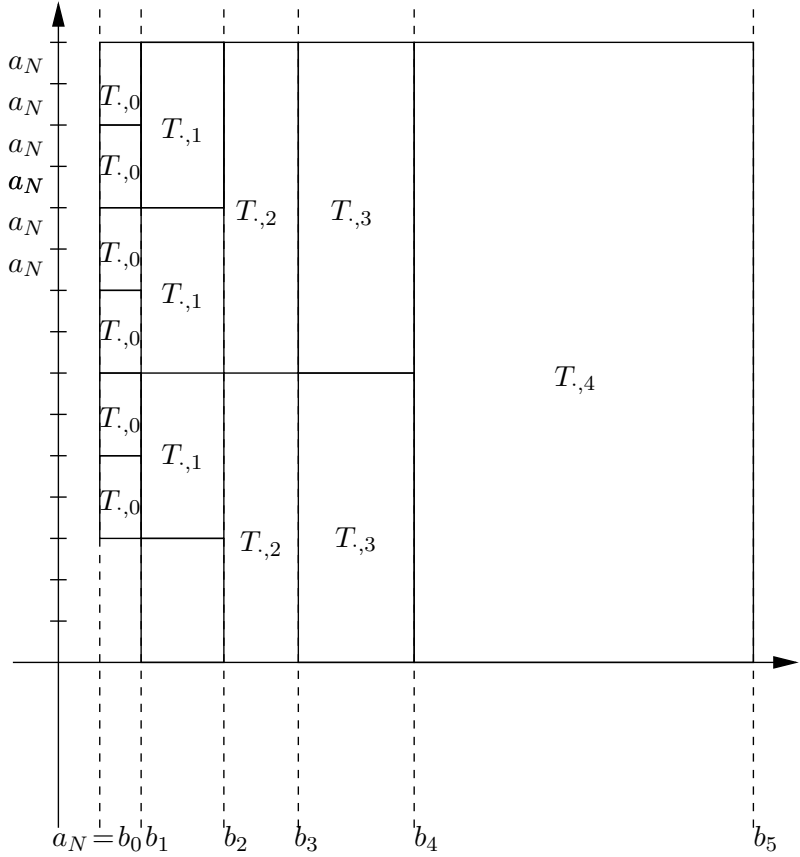}
    \caption{The tiles $T_{I,j}$}
    \label{fig:tiles}
\end{figure}

  The family $\FF$ gives rise to a  family $\TT$ of right halves of the Carleson rectangles  in $\FF$, which we will call {\em tiles}, and which will form a decomposition of the closed half plane $\overline{\CC_{a_N}}$ into disjoint sets.
   These are rectangles of the form $T_{I,j}=  [b_{j}, b_{j+1}) \times I$, where $Q_{I,j} \in \FF$, see Figure \ref{fig:tiles}. 

We say that an interval $I \in \DD_N$ belongs to generation $j$, if it is the base for a Carleson square $Q_{I,j} \in \FF$. Thus each $I \in \DD_N$ belongs to at least one, and possibly more than one,
generation.

The idea of the construction below is to define $\mu_N$ as the ``largest possible part" of $\mu$ which can be dominated by $\nu_N$, in terms of Condition (\ref{cond:ncarl}).
Recall that by the  Carleson-type Condition (\ref{eq:zencond})  we are given, we have in particular
$$
   \mu(T_{I,0})=   \mu(Q_{I, 0}) \le \nu(Q_{I, 0}) =  \nu_N(Q_{I, 0}) 
$$
for all intervals $I$ of generation $0$ in  $\DD_N$.
For such intervals, we define the  remaining part of $\nu_N$ by
$$
        \nu^0_N|_{ \{a_N\} \times I   }  =   \frac{\nu_N(Q_{I,0}) - \mu(T_{I, 0})}{\nu_N(Q_{I,0}) } \nu_N |_{ \{a_N\} \times I   } .
$$
This defines a measure $\nu_N^0$ on $ \{a_N\} \times i \RR$.

 In the second step, we define a measure $\nu_N^1$ on   $ \{a_N\} \times i \RR$ by letting
 $$
            \nu^1_N|_{ \{a_N\} \times I'   }  =  \left\{     \begin{matrix}     
 \frac{\nu^0_N(Q_{I',1}) - \mu(T_{I', 1})}{\nu^0_N(Q_{I',1}) } \nu_N |_{ \{a_N\} \times I'   } & \text{ if } &  \nu^0_N(Q_{I',1})  >  \mu(T_{I', 1}) ,  \\  
                                                                                                                                                                  0& \text{ if } &  \nu^0_N(Q_{I',1})  \le  \mu(T_{I', 1})  , \\
            \end{matrix}     \right.
 $$
for intervals $I' \in \DD_N$ of generation $1$. In the next and all following steps, having already defined the 
measure $\nu_N^j$ on $ \{a_N\} \times i \RR$ for some $j \ge 0$, we let
 \begin{multline}   \label{eq:measdef}
            \nu^{j+1}_N|_{ \{a_N\} \times J  }  = \\
               \left\{     \begin{matrix}      \frac{\nu^j_N(Q_{J,j+1}) - \mu(T_{J, j+1})}{\nu^j_N(Q_{J,j+1}) } \nu_N |_{ \{a_N\} \times J  } & \text{ if } &  \nu^j_N(Q_{J,j+1})  >  \mu(T_{J, j+1}) ,  \\  
                                                                                                                                                                  0& \text{ if } &  \nu^j_N(Q_{J,j+1})  \le  \mu(T_{J, j+1}),  \,\, \nu^j_N(Q_{J,j+1}) >0 ,  \\
                                                                                                                                                                   0& \text{ if } &    \nu^j_N|_{ \{a_N\} \times J  }     =0 ,  \\
  \end{matrix}            \right.
 \end{multline}
for intervals of generation $j+1$ in $\DD_N$, thereby defining the measure $\nu^{j+1}_N$ on $ \{a_N\} \times i \RR$. Depending on which case in (\ref{eq:measdef}) appears, we say that $(J,j+1)$ is of type 1, 2 or 3,
respectively. 
We are finally ready to define the measure $\mu_N$. We start by letting
$$
      \mu_N|_{T_{I,0}}   =   \mu|_{T_{I,0}}
$$
for intervals $I \in \DD_N$ of generation $0$. For $T_{I, j+1} \in \TT$, $j \ge 0$, let
\begin{multline*}
       \mu_N|_{T_{I,j+1}}         = \\
        \left\{     \begin{matrix}      \mu |_{ T_{I,j+1} } & \text{ if } &  \nu^j_N(Q_{I,j+1})  >  \mu(T_{I, j+1}) ,  \\  
                                                                                                                                                                     \frac{ \nu_N^{j}(Q_{I, j+1})}{\mu(T_{I, j+1})}  \mu |_{ T_{I,j+1}  } & \text{ if } &  \nu^j_N(Q_{I,j+1})  \le  \mu(T_{I, j+1}),  \nu^j_N(Q_{I,j+1}) >0 ,  \\
                                                                                                                                                                   0& \text{ if } &    \nu^j_N|_{ T_{I,j+1} }     =0.   \\
                                                                                                                                                                        \end{matrix}     \right.
 \end{multline*}                                                                                                                                                                
   Since the tiles in $\TT$ form a decomposition of $\overline{\CC_{a_N}}$, this defines a measure $\mu_N$ on    $\overline{\CC_{a_N}}$.
   
 By construction, we have  for each $ Q_{I,j} \in \FF$, $(I,j)$ of type 2 or type 3:
    $$
           \mu_N(Q_{I, j})   =   \nu_N (\{a_N\} \times I).
    $$       
    If $(I,j)$ is of type 1, then
    $$
        \mu_N(Q_{I, j})   <   \nu_N (\{a_N\} \times I).
    $$                                                                                                                              
Hence $\mu_N$ satisfies Condition 3 in Theorem \ref{thm:Zen} with respect to $\nu_N$. Let us now look at the Carleson condition for $\mu - \mu_N$.

If $(I,j)$ is of type 1, then $\mu_N(T_{I,j}) = \mu(T_{I,j})$ by construction, and $(\mu - \mu_N)(T_{I,j})=0$. Therefore, by decomposing the Carleson square
$Q_{I,j}$ into 
$$
    Q_{I,j} = T_{I,j}  \cup \left( \bigcup_{I' \subset I, Q_{I', j-1} \in \FF}    Q_{I', j-1}    \right)
$$
and iterating if necessary, we see that we have to check Condition 3 only for Carleson squares $Q_{I,j}$ with $(I,j)$ of type 2 or type 3. But in this case,
  $$
           \mu_N(Q_{I, j})   =   \nu_N (\{a_N\} \times I)
    $$     
and thus  by the original Condition 3 in Theorem \ref{thm:Zen} for $\mu$ and $\nu$,
$\mu- \mu_N$ satisfies Condition 3 with respect to $\nu - \nu_N$, again with constant $1$ on Carleson squares in $\FF$. 

Finally, again by  Condition (\ref{eq:zencond}) in Theorem \ref{thm:Zen} for $\mu$ and $\nu$, we see that for all tiles  $T_{I,j}$ contained in the strip $\{ z \in \CC: a_N \le \Re z < a_{N+1}\}$, $ \mu_N|_{T_{I,j}} = \mu|_{T_{I,j}}$. Therefore, $\mu- \mu_N$ is supported on the closed half-plane $\overline{ \CC_{a_{N+1}}}$. We can now make an induction step by applying the same procedure to the measures
$\mu- \mu_N$, $\nu- \nu_N$ with respect to the half-plane $\overline{ \CC_{a_{N+1}}}$ to construct $\mu_{N+1}$, etc. We thus obtain a decomposition
$$
        \mu = \sum_{n=N}^\infty \mu_n
$$
satisfying the Conditions (\ref{cond:sum}), (\ref{cond:supp}) and (\ref{cond:ncarl}) in the lemma.

It remains to be shown that for each $n > N$,
$\mu_n$ is a Carleson measure with respect to the shifted half-plane  $\CC_{a_{n-1}}$, with the appropriate estimate of the Carleson constant. Let $Q$ be a Carleson square in $\CC_{a_{n-1}}$
 over the interval $I$. Recall that by (\ref{cond:ncarl}), there exists $C' >0$ with $\mu_n(\tilde Q) \le C' \nu_n(\tilde Q)$ for each Carleson square $\tilde Q $ in $\CC_+$ and that
 $(a_n - a_{n-1}) \ge c a_n$.
 
 If the sidelength $|I|$ of $Q$ is less than $a_n - a_{n-1}$, then $Q$ has empty intersection with the support of $\mu_n$, and 
$\mu_n(Q) =0$. If $|I| \ge   a_n - a_{n-1}   $, then $Q$ can be covered by a Carleson square $\tilde Q$ in $\CC_+$ with sidelength at most $\frac{1}{c} |I|$. Thus
$$
    \mu_n(Q) \le \mu_n(\tilde Q) \le C' \nu_n(\tilde Q) \le C' \frac{1}{c}  \beta_n  | I |  ,
$$
and we obtain the desired result. The result for $\mu_N$ is shown in the same way.
\qed

Now let $f \in A_\nu^p$. Note that by the $(\Delta_2)$-condition, each of the norms $\|f\|_{H^p_{a_n}}$, denoting the norms on the Hardy spaces $H^p(\CC_{a_n})$ of the shifted half planes $ \CC_{a_n}$, is finite. 
Restricting $\mu$ to some closed half-plane $\overline{\CC_{a_N}}$, and using the decomposition in Lemma \ref{lemm:mainest}, we obtain
\begin{eqnarray*}
&&    \int_{\overline{\CC_{a_N}}} |f(z)|^p d\mu(z)  \\
  &= & \sum_{n =N}^\infty      \int_{\CC_+}   |f(z)|^p d\mu_n(z)    \\
  &\le & C_p \frac{C'}{c}     \sum_{n = N}^\infty     \beta_n    \|f\|_{H^p_{a_{n-1}}}^p   \\
  &\le & (2 R)^3   C_p \frac{C'}{c}     \sum_{n = N}^\infty     \beta_{n-1}    \|f\|_{H^p_{a_{n-1}}}^p   \\ 
  & =  &    (2 R)^3   C_p \frac{C'}{c}      \sum_{n = N}^\infty    \tilde \nu([a_{n-1}, a_{n}))  \|f\|^p_{H^p_{a_{n-1}}} \\
    & \le &    (2 R)^3   C_p \frac{C'}{c}     \|   f(z)\|^p_{A_\nu^p}.         \\
\end{eqnarray*}
Here, we use that for any $f \in A^p_{\nu}$, the map
$$
       r \mapsto \int_{-\infty}^\infty   | f(r + it)|^p dt
$$
is non-increasing. Using that
$$
           \int_{\CC_+} |f(z)|^p d\mu(z) = \sup_{N \in \ZZ}    \int_{\overline{\CC_{a_N}}} |f(z)|^p d\mu(z),
$$
we obtain the desired estimate.

Now we can finish the proof of the theorem by showing that a sequence of positive numbers  $(a_n)_{n \in \ZZ}$ with the required properties (\ref{eq:delta2}), (\ref{eq:delta2s})  exists, and that we can also treat the case $\nu(\{0\}) >0$.
 Let $R$ be the $(\Delta_2)$ constant of the measure $\tilde \nu$, and 
let $F$ be the function given by,
$$
   F: [0, \infty) \rightarrow \RR, \quad F(r) = \tilde \nu([0, r)).
$$
$F$ is left continuous, and also right continuous up to countably many jumps. By the $(\Delta_2)$-condition, 
$$
         F(2t) \le R F(t), 
$$
so in particular, a jump of $F$ at $t$ may be no more than $(R-1)F(t)$.

If $\tilde \nu(\{0\}) = 0$, then let for $n \in \ZZ$
$$
         a_n = \sup\{r \ge 0: F(r) \le (2R)^{2n}  \},
$$
provided that the supremum is finite. If the supremum is infinite, we stop the sequence at the corresponding $n$.

By the right continuity of $F$, $a_n >0$ for all $n \in \ZZ$, and by the condition on jumps of $F$,
$$
       (2R)^{2n} \ge       F(a_n)  \ge \frac{1}{R}  (2R)^{2n}
$$
and therefore
$$
   \frac{\tilde \nu([0, a_{n+1})) }{\tilde \nu([0, a_{n})) }  \ge \frac{ (2R)^{2n+2}}{R (2R)^{2n}}  \ge 4R.
$$
Hence $a_{n+1} \ge 2 a_n$ and
$$
   \frac{a_{n+1} - a_n}{a_{n+1}} \ge \frac{1}{2}.
$$
Furthermore,
$$
                    (2R)^{2n}  (4R -1 )     \le          \tilde \nu([a_n, a_{n+1}))   =  F(a_{n+1}) - F(a_n) \le   (2R)^{2n +2} ,
$$
hence
$$
                   2R \le       \frac{\tilde \nu([a_n, a_{n+1})) }{ \tilde \nu([a_{n-1}, a_{n})) }   \le    (2R)^3.
$$

If  $\tilde \nu(\{0\}) \neq 0$, then let $a_0=0$, $\beta_0 = \tilde \nu(\{0\})  $, and let 
$$
         a_n = \sup\{r \ge 0: F(r)  \le (2R)^{2(n+1)}   \tilde \nu(\{0\})    \}    \text{ for } n \in \NN   .
$$
We see in the same way as before that properties (\ref{eq:delta2}), (\ref{eq:delta2s}) hold, if we replace $\tilde \nu([a_0, a_1))$ by $ \tilde \nu((a_0, a_1)) $. We write
$\beta_1 = \tilde \nu((a_0, a_1))$ and $\beta_n = \tilde \nu([a_{n-1}, a_{n}))$ for $n \ge 2$. Then
\begin{multline*}
    \beta_0   \|f\|^p_{H^p_{a_{0}}}  +     \sum_{n =1}^\infty    \beta_n   \|f\|^p_{H^p_{a_{n+1}}}    \le  \int_{\CC_+} |f(z)|^p d\nu(z)   \le    \beta_0   \|f\|^p_{H^p_{a_{0}}} +   \sum_{n =1}^\infty \beta_n   \|f\|^p_{H^p_{a_{n}}},
\end{multline*}
and the same construction as before applies. If the sequence $(a_n)$ is finite to the right, an analogous argument can be made.
\qed


The following proposition is elementary and appears for special cases in \cite{zen09,zen10}. Partial results are also given in
\cite{DP94,DGM}.
\begin{proposition} Let $A^2_\nu$ be a Zen space,   and let $w:(0, \infty) \rightarrow \RR_+$ be given by
$$
    w(t) = 2 \pi \int_0^\infty  e^{-2rt} d \tilde \nu(r)   \qquad (t >0).
$$
Then the Laplace transform defines an isometric map $\LL: L^2_w(0, \infty) \rightarrow A^2_\nu$.
\end{proposition}
Note that the existence of the integral is guaranteed by the $(\Delta_2)$-condition.
\bigskip

\proof Let $f \in L^2_w(0, \infty)$. Then
\begin{eqnarray*}
&&\sup_{\eps >0} \int_{\CC_+} |\LL f(z+\eps)|^2 d\nu(z)\\
& =& 
\sup_{\eps >0}  \int_0^\infty  \| (\LL f)(\eps + r + \cdot)\|^2_{L^2(i \RR)} d \tilde \nu(r)\\
&=& \sup_{\eps >0}  \int_0^\infty  \|\FF (e^{-(r+\eps)\cdot} f)\|^2_{L^2(\RR)} d \tilde \nu(r)\\
&=& \sup_{\eps >0}  \int_0^\infty  2 \pi \| e^{-(r+\eps)\cdot} f\|^2_{L^2(0, \infty)} d \tilde \nu(r)\\
&=& \sup_{\eps >0}  \int_0^\infty |f(t)|^2 2 \pi \int_0^\infty e^{-2(r+\eps)t} d\tilde \nu(r) dt \\
&=& \int_0^\infty |f(t)|^2 w(t) dt
\end{eqnarray*}
by isometry of the Fourier transform and the dominated convergence theorem.
\qed

Here is a Laplace--Carleson Embedding Theorem, which is an immediate consequence.

\begin{theorem} \label{thm:CarlLap}          Let $A^2_\nu$ be a Zen space, $ \nu = \tilde \nu \otimes d \lambda$, and let $w:(0, \infty) \rightarrow \RR_+$ be given by
\begin{equation}\label{massw}
    w(t) = 2 \pi \int_0^\infty  e^{-2rt} d \tilde \nu(r)   \qquad (t >0).
\end{equation}
Then the following are equivalent:
\begin{enumerate}
\item The Laplace transform $ \LL$ given by  $ \LL f (z) = \int_0^\infty e^{-t z} f(t) dt$ defines a bounded linear map
$$
\LL: L^2_{w}(0,\infty) \rightarrow L^2(\CC_+, \mu).
$$
\item
For a sufficiently large $N \in \NN$, there exists a constant $C >0$ such that
$$
      \int_{\CC_+}    \left|  (\LL t^{N-1} e^{-\lambda t }) (z)\right|^2 d \mu(z)  \le C     \int_0^\infty    | t^{N-1} e^{-\lambda t } |^2  w(t) dt         \text{ for each } \lambda \in \CC_+.
$$
 \item There exists a constant $C>0$ such that 
\begin{equation*}   
 \mu(Q_I) \le C \nu(Q_I) \text{ for each Carleson square }Q_I.
 \end{equation*}
\end{enumerate}
\end{theorem}
\proof Noticing that $\LL  (t^{N-1} e^{-\lambda t } )$ is a scalar multiple of $(k_\lambda)^N$, this follows immediately from Theorem \ref{thm:Zen}.   \qed


\subsection{Hankel operators on Zen spaces}
\label{sec:2.2}

The boundedness of little Hankel operators on weighted Bergman spaces has been studied in various settings, see e.g. \cite{DP94}, \cite{zhu1}, \cite{janson}, \cite{yong},
though mostly for standard weights. Generally, boundedness results for Hankel operators are
often  connected to Carleson embedding results, at least concerning sufficient conditions for boundedness.

We will show in this subsection that the boundedness of Hankel operators on Zen spaces can be deduced from a Carleson
measure condition on the symbol. In the case of the Hardy space or the standard-weighted Bergman spaces, this condition is also necessary.

\begin{theorem}     \label{thm:hankel}
Let $A^2_\nu$ be a Zen space, $d\nu = d\tilde \nu \otimes d\lambda$.
Let $b: \CC_+ \rightarrow \CC$ be analytic, $b \in H^2(\CC_+)$. 
If the measure 
\begin{equation}    \label{eq:hankel}
    |b'(z)|^2  \Re z   \: \tilde \nu([0, \Re z) )   \, dA(z)
\end{equation}
is a $\nu$-Carleson measure on $\CC_+$, then the little Hankel operator
$$ A^2_\nu \rightarrow \overline{A^2_\nu},  \quad f \mapsto Q_\nu \bar b f  
$$
is bounded.
Here, $Q_\nu$ denotes the orthogonal projection $Q_\nu: L^2(\CC_+, d \nu) \rightarrow \overline{A^2_\nu}$.

\end{theorem}

\vspace{0.5cm}
\proof
For the proof,  we follow the lines of
 the proof of the Fefferman-Stein Duality Theorem in \cite{garnett}.
Suppose that $| b'(z)|^2  \Re z \,   \tilde \nu([0, \Re z ))    dA(z)$ is a $\nu$-Carleson
measure and let $f,g \in A^2_\nu$. Then
\begin{eqnarray*}
\langle Q_\nu \bar b f, \bar g \rangle_{ A^2_\nu} &=& 
\int_0^\infty  \int_{i \RR} \bar b(r+t) f(r+t) g(r+t) dt d \tilde \nu(r) \\
&=& \int_0^\infty \int_{\CC_+} \bar b'(z+r)    ( f(r+ \cdot) g(r+ \cdot))'(z) \Re z dA(z) d \tilde \nu(r),\\
\end{eqnarray*}
where we use the Littlewood--Paley identity
$$
  \int_0^\infty \int_{i \RR}|f(r+t)|^2 dt d \tilde \nu(r)  = \int_0^\infty \int_{\CC_+} 
     |f'(z+r)|^2  \Re z dA(z) d \tilde \nu(r)
$$
and its polarization. Now
\begin{eqnarray*}
\lefteqn{\left|\int_0^\infty \int_{\CC_+} \bar b'(z+r)     f(r+ z) g'(r+ z) \Re z dA(z) d \tilde \nu(r)\right|}\\
& \le& (\int_0^\infty \int_{\CC_+} |\bar b'(z+r)|^2 |f(r+ z)|^2 \Re z dA(z)d \tilde \nu(r))^{1/2}\\
&& \qquad \qquad\qquad
       (\int_0^\infty \int_{\CC_+} |g'(z+r)|^2 \Re z dA(z)d \tilde \nu(r))^{1/2} \\
&=& ( \int_{\CC_+} |\bar b'(z)|^2 |f(z)|^2 [\int_0^{\Re z} \Re (z-r) d \tilde \nu(r)]dA(z))^{1/2}\\
&& \qquad\qquad\qquad
(\int_0^\infty \int_{\CC_+} |g'(z+r)|^2 \Re z dA(z)d \tilde \nu(r))^{1/2} \\
&\lesssim & \|f\|_{A^2_\nu} \|g\|_{A^2_\nu},
\end{eqnarray*}
since
$$
|b'(z)|^2  \int_0^{\Re z} \Re (z-r) d \tilde \nu(r) dA(z) \le   |b'(z)|^2  \Re z \: \tilde \nu([0, \Re z) )  \, dA(z)
$$
and the latter is a $\nu$-Carleson measure.
The second term \\ $\int_0^\infty \int_{\CC_+} \bar b'(z+r)f'(r+ z) g(r+ z) \Re z dA(z) d \tilde \nu(r)$
is estimated accordingly. Hence the Hankel operator is bounded.
\qed

In the case of the Hardy space $H^2(\CC_+)$, that is, if $\tilde \nu$ is the Dirac measure in $0$, the Carleson condition on the measure
(\ref{eq:hankel}) is a well-known characterisation of the boundedness of the Hankel operator, see e.g. \cite{garnett}.

For further applications of the theorem above, recall the Bloch space $\BB(\CC_+)$ on $\CC_+$,
$$
 \BB(\CC_+) = \left\{   f \in Hol(\CC_+):       \| f \|_{\BB}=  \sup_{z \in \CC_+} |f'(z)| \Re z < \infty       \right\}.
$$
Let us consider the following inverse $(\Delta_2)$ condition on $\tilde \nu$:
\begin{equation}   \label{eq:invdelta}
     \sup_{M >1}   \,   \inf_{r >0}   \,  \frac{\tilde \nu([0, Mr )) }{ \tilde \nu([0, r ))   }   >1.
\end{equation}
\begin{corollary}    \label{cor:hankel}
Let $A_\nu^2(\CC_+)$ be a Zen space with $\tilde \nu$  satisfying the inverse $(\Delta_2)$ condition (\ref{eq:invdelta}) and let $b \in \BB(\CC_+)$.
Then the little Hankel operator
$$ A^2_\nu \rightarrow \overline{A^2_\nu},  \quad f \mapsto Q_\nu \bar b f  
$$
is bounded.
\end{corollary}
\proof Let $b \in \BB(\CC_+)$ and let $Q_I$ be a Carleson square. We write $F(r) = \tilde \nu ([0,r))$ for $r \ge 0$. Then 
\begin{eqnarray*}    
   \int_{Q_I}       |b'(z)|^2    \Re z  \,   \tilde \nu( [0, \Re z) ) \, dA(z) & \le  & \|b\|^2_{\BB}     \int_{Q_I}     \frac{1}{\Re z} \int_0^{\Re z}  d \tilde \nu(r) dA(z) \\
   &\le & |I|    \int_{0}^{|I|}   \frac{1}{s} F(s) ds. \\
\end{eqnarray*}

To finish the proof, it remains only to be shown that there exists $C>0$ with
$$
 \int_{0}^{x}   \frac{1}{s} F(s) ds   \le C  F(x) \text{ for all } x \in [0, \infty).
$$
Choosing $M , \gamma > 1$ with 
$$
    \inf_{r >0}    \frac{\tilde \nu([0, Mr )) }{ \tilde \nu([0, r ))   }   \ge \gamma,
$$
we find that
\begin{multline*}
         \int_{0}^{x}   \frac{1}{s} F(s) ds  \le \sum_{n=0}^\infty        \int_{M^{-(n+1)} x}^{M^{-n} x}   \frac{1}{s} F(s) ds \\
         \le  \sum_{n=0}^\infty \frac{1}{x}  (M^{-n} - M^{-(n+1)} )   x M^{n+1}   F(M^{-n} x)    \\ \le (M-1)  \sum_{n=0}^\infty  \gamma^{-n}  F(x)    = \frac{\gamma (M-1) }{\gamma-1} F(x).
\end{multline*}
Hence the $\nu$-Carleson condition for the measure (\ref{eq:hankel}) in Theorem \ref{thm:hankel} applies, and the little Hankel operator is bounded.
\qed

{\bf Remark.}  In the case of the standard weighted Bergman spaces $A^2_\alpha(\CC_+)$ with $\alpha >-1$,   we have $d \tilde \nu(r) =  r^\alpha dr$, and the inverse $(\Delta_2)$ condition (\ref{eq:invdelta}) holds.
 It is well-known that 
the little Hankel operator with symbol $b$ is bounded on $A^2_\alpha$, if and only if $b \in \BB(\CC_+)$ (see e.g. \cite{zhu}). Therefore, the $\nu$-Carleson condition
 for the measure (\ref{eq:hankel}) in Theorem \ref{thm:hankel}, respectively the Bloch condition on the symbol in Corollary (\ref{cor:hankel}),  are also necessary in this case.

\section{$L^p-L^q$ embeddings}    \label{sec:pq}
As mentioned above, there is no known full characterization of boundedness
of Laplace--Carleson embeddings 
$$
  L^p(0,\infty) \rightarrow L^{q}(\CC_+, \mu), \quad f
  \mapsto \LL f= \int_0^\infty e^{-t \cdot} f(t) dt.
$$
However, characterizations are possible in case $p\le 2$, if the conjugate index $p'$ satisfies $p' <q$, and in some cases with additional information on the support on the measure.
We list some results for natural spectral inclusion conditions which appear naturally which appear naturally in the theory of operator
semigroups.
In the cases we consider here, the oscillatory part of the Laplace transform can be discounted,
and a full characterization of boundedness can be achieved.

First, let us make the following simple observation.

\begin{proposition}   \label{prop}    Let $\mu$ be a positive regular Borel measure on $\CC_+$, let $1 \le p,q < \infty$ and suppose that
the Laplace--Carleson embedding
$$
     \LL: L^p(0,\infty) \rightarrow L^q(\CC_+, \mu), \quad f \mapsto \LL f,
$$
is well-defined and bounded.
Then there exists a constant $C_{p,q}>0$ such that for all intervals $ I \subset i\RR$,
\begin{equation}  \label{eq:duren}
  \mu(Q_I)  \le  C_{p,q}  |I|^{q/p'}    \text{ if } p>1,  \qquad
   \mu(Q_I)   \le  C_{p,q}  \text{ if } p=1.
\end{equation}

\end{proposition}
\proof Let $Q_I$ be a Carleson square with centre $\lambda_I$. Note that
$$
   \|(\LL e^{- \cdot \bar \lambda_I})\|^q_{L^q(\CC_+, \mu)} \ge \int_{Q_I} |(\LL e^{- \cdot \bar \lambda_I})|^q d \mu
      \ge \frac{1}{(4 \re \lambda_I)^q} \mu(Q_I) = \frac{1}{2^q |I|^q} \mu(Q_I),
$$ 
and
$$
   \|e^{- \cdot \lambda_I}\|^p_p = \frac{1}{p \re \lambda_I} = \frac{2}{p |I|}
   $$
for $1\le p,q < \infty$.
Hence there exists a constant $C_{p,q} >0$ such that
$$
        \mu(Q_I)    \le C_{p,q} {|I|^{q/p'}} \text{ if } p >1 \text{ and }   \mu(Q_I)    \le C_{p,q}  \text{ if } p =1.
$$
This concludes the proof. \qed

The proposition immediately yields the following theorem:

\begin{theorem}
 \label{thm:pqexp} Let $\mu$ be a positive regular Borel measure supported in the right half-plane 
 $\CC_{+} $, and let $1 < p' \le q < \infty$, $p\le2$.
 Then the following are equivalent:
\begin{enumerate}
\item The Laplace--Carleson embedding
$$
     \LL: L^p(0,\infty) \rightarrow L^q(\CC_+, \mu), \quad f \mapsto \LL f,
$$
is well-defined and bounded.
\item There exists a constant $C>0$ such that 
\begin{equation}  
\mu(Q_I) \le C
  |I|^{q/p'} \text{ for all intervals }   I \subset i\RR
  \end{equation}
   \item  There exists a constant $C>0$ such that $\|\LL e^{- \cdot z}
  \|_{L^q_\mu} \le C \|e^{- \cdot z}\|_{L^p}$ for all $z \in \CC_{+}$.
  \end{enumerate}
\end{theorem}
\proof
Obviously (1) $\Rightarrow$ (3), and (3) $\Rightarrow$ (2) by the Proposition.

 This leaves (2) $\Rightarrow$ (1), which is also easy:
By the Hausdorff--Young inequality, the map $\LL: L^p(0, \infty) \rightarrow L^{p'}(i \RR)$ is bounded. By Duren's theorem, the Poisson extension
$L^{p'}(i \RR) \rightarrow L^q_\mu(\CC_+)$ is bounded, given the Carleson condition (\ref{eq:duren}). The composition of both gives the boundedness of the Laplace transform
$
   \LL: L^p(0,\infty) \rightarrow L^q(\CC_+, \mu).
$
This concludes the proof.\qed

The case $p >2 $ is much more subtle: in particular the Laplace transform no longer maps $L^p(0,\infty)$ into $H^{p'}(\CC_+)$.
We give two special cases here, that of the measure $\mu$ being supported in a strip and that of $\mu$ being supported in a sector.


\subsection{Sectorial measures}

 If the measure $\mu$ is supported on a sector $S(\theta) = \{
z \in \CC_+: |\arg z| < \theta \}$ for some $0< \theta <
\frac{\pi}{2}$,
then the oscillatory part of the Laplace transform can be discounted,
and a full characterization of boundedness can be achieved (see also
\cite{haak}, Theorem 3.2 for an alternative characterization by means
of a different measure).

\begin{theorem} \label{thm:pqemb} Let $\mu$ be a positive regular Borel measure supported in a
  sector $S(\theta) \subset \CC_+$, $0 < \theta < \frac{\pi}{2}$, and
  let $q \ge p >1$. Then the following are equivalent:
\begin{enumerate}
\item The Laplace--Carleson embedding
$$
     \LL: L^p(0,\infty) \rightarrow L^q(\CC_+, \mu), \quad f \mapsto \LL f,
$$
is well-defined and bounded.

\item There exists a constant $C>0$ such that $\mu(Q_I) \le C
  |I|^{q/p'}$ for all intervals in $I \subset i\RR$ which are symmetric
  about $0$.

\item  There exists a constant $C>0$ such that $\|\LL e^{- \cdot z}
  \|_{L^q_\mu} \le C \|e^{- \cdot z}\|_{L^p}$ for all $z \in \RR_+$.
  
  \item  There exists a constant $C>0$ such that $\|\LL e^{- \cdot 2^n}
  \|_{L^q_\mu} \le C \|e^{- \cdot 2^n}\|_{L^p}$ for all $n \in \NN$.

\end{enumerate}

\end{theorem} 
\proof (2) $\Rightarrow$ (1)
For $n \in \ZZ$, let 
\[
T_n = \{ x+ iy \in \CC_+: 2^{n-1} < x \le 2^{n}, -2^{n-1}< y  \le 2^{n-1} \}.
\]
  That is, $T_n$ is the right half of the Carleson square $Q_{I_n}$
  over the interval $I_n = \{ y \in \RR, |y| \le 2^{n-1} \}$. The
  $T_n$ are obviously pairwise disjoint.

Without loss of generality we assume  $0<\theta < \arctan(\frac{1}{2}) $, in which
case $S(\theta) \subseteq \bigcup_{n = - \infty}^\infty T_n$.

Now let $z \in T_n$ for some $n \in \ZZ$. Then we obtain, for $f\in L^p(0,\infty)$,
$$
 |\LL f(z)| \le \int_0^\infty |e^{-zt}| |f(t)| dt \le \int_0^\infty
 |e^{-2^{n-1}t}| |f(t)| dt \le C_\Theta 2^{-n+1} Mf(2^{-n+1}),
$$
where $C_\Theta >0$ is a constant dependent only on the integration kernel $\Theta(t)= \chi_{[0,\infty)}(t+1) e^{-t-1}$ and 
$Mf$ is the Hardy--Littlewood maximal function. We refer to e.g. \cite{stein}, page 57, equation (16) for a pointwise estimate between the maximal
function induced by the kernel $\Theta$ and $M$. We can easily dominate $\Theta$ by a positive, radial, decreasing $L^1$ function here. 
Consequently,
\begin{eqnarray*}
     \int_{S(\theta)} |\LL f(z)|^q d\mu(z)  
&\le& \sum_{n= -
     \infty}^\infty 2^{q(-n+1)} (Mf(2^{-n+1}))^q \mu(T_n)\\
& \le &  C_\Theta^q \sum_{n= -
     \infty}^\infty 2^{q(-n+1)} 2^{n q/p'} (Mf(2^{-n+1}))^q  \\
 &  = &   C_\Theta^q\sum_{n= -
     \infty}^\infty 2^{q/p'} (2^{(-n+1)} (Mf(2^{-n+1})^p))^{q/p} \\
& \le & C_\Theta^q  2^{q/p'}  \left(\sum_{n= -
     \infty}^\infty  2^{(-n+1)} (Mf(2^{-n+1}))^p\right)^{q/p} \\
& \lesssim & \|f\|_{L^p}^{q}.
\end{eqnarray*}
Note that in the case $1 < p \le 2$, (but not for $p>2$) this result can also easily be
deduced from the Hausdorff--Young inequality and Duren's Theorem \cite{duren}.\\

(4) $\Rightarrow$ (2) Let $I \subset i\RR$ be an interval which is symmetric about $0$. We can assume without loss that $|I| = 2^n$. It is easy to see that
$$
   |(\LL e^{- 2^{n-1}\cdot}) (z)| = \left|\frac{1}{2^{n-1} + z}\right| \ge \frac{1}{2^{n+1}} \text{ for } z \in Q_I.
$$ 
Thus
\begin{eqnarray*}
 \mu(Q_I) &\le& (2^{n+1})^q \int_{\CC_+} | (\LL e^{-2^{n-1})(z) \cdot}|^q d \mu(z) \\
&\le& C (2^{n+1})^q \|e^{-2^{n-1} \cdot}\|^{q}_p
  \approx 2^{nq} 2^{-nq/p} = 2^{nq/p'}.
\end{eqnarray*}
(1) $\Rightarrow$ (3)  and (3) $ \Rightarrow$ (4) are obvious.
\qed

\begin{remark}
Let $\mu$, $\theta$, $p$ and $q$ be as in Theorem \ref{thm:pqemb}. 
In \cite{haak}, Theorem 3.2, essentially the equivalence of the following statements is shown for discrete measures:
\begin{enumerate}
\item The Laplace--Carleson embedding
$$
     \LL: L^p(\RR_+) \rightarrow L^q(\CC_+, \mu), \quad f \mapsto \LL f,
$$
is well-defined and bounded.
\item There exists a constant $C>0$ such that $\tilde \mu(Q_I) \le C
  |I|^{q/p}$ for all intervals in $I \subset i\RR$ which are symmetric
  about $0$, where $d\tilde \mu(z) = |z|^q d\mu(\frac{1}{z})$.
  \end{enumerate}
The proof also uses the method of the maximal function.
\end{remark}


Now let us consider the case $p>q$ for sectorial measures $\mu$. We will, among others, obtain a condition in terms of the balayage
$S_{ \mu}$ of $ \mu$
(compare this with the characterization of bounded $H^p(\CC_+) \rightarrow L^q(\CC_+, \mu)$
embeddings for $p>q$ in \cite{l2}). Recall that the balayage $S_\mu$ of a positive Borel measure $\mu$ on $\CC_+$ is given by
$S_\mu(t)=\int_{\CC_+} p_z(t) d \mu(z)$, where $p_z$ is  the Poisson kernel, defined in the introduction  as
\[
p_z(t)=\frac{1}{\pi} \frac{y}{x^2+(y-t)^2}, \qquad (z=x+iy \in \CC_+, \  t \in \RR).
\]

To look at a dyadic version, let
\[
T_n 
= \{ x+ iy \in \CC_+: 2^{n-1} < x \le 2^{n}, -2^{n-1}< y  \le 2^{n-1}
\}
\] as in the previous proof, and for $k \in \ZZ$, let $T_{n,k} = T_n+ i k
2^n$, so the $T_{n,k}$ are the translates of $T_n$ parallel to the
imaginary axis. Similarly, let
$I_{n,k} = I_n+ k 2^n$ be the translates of $I_n:=\{y\in \mathbb R: |y|\le 2^{n-1}\}$. The $\{
T_{n,k}:n,k \in \ZZ \}$ then form a dyadic tiling of the right
half plane. We write $S_n = \cup_{k \in \ZZ} T_{n,k} = \{ z \in \CC: 2^{n-1} < \Re z \le 2^n \}$.  Let
$$
   S^d_{\mu}(t) = \sum_{n,k \in \ZZ} \chi_{I_{n,k}}(t) \frac{ \mu(T_{n,k})}{2^n}.
$$
$S^d_\mu$ is called the dyadic balayage of $\mu$. Note that $S^d_{ \mu} \le 2 \pi S_{ \mu}$ pointwise, since
\begin{multline}   \label{eq:upperbalay}
    S_{ \mu}(t) = \int_{\CC_+} p_z(t) d  \mu(z) \ge
    \sum_{n=-\infty}^\infty  \sum_{k=-\infty}^\infty \chi_{I_{n,k}}(t)
    \int_{T_{n,k}} p_z(t) d  \mu(z)  \\
\ge \sum_{n=-\infty}^\infty  \sum_{k=-\infty}^\infty \chi_{I_{n,k}}(t)
    \inf_{z \in T_{n,k}}\{p_z(t)\}  \mu(T_{n,k}) \ge \frac{1}{2 \pi} \sum_{n,k \in \ZZ} \chi_{I_{n,k}}(t) \frac{ \mu(T_{n,k})}{2^n}.
\end{multline}

In the special case that $\mu$ is sectorial with opening angle $\theta < \pi/2$, the measure $\mu$ is supported on $\cup_{n \in \ZZ} T_n$, and we get a particularly simple
form of the dyadic balayage $S^d_\mu$, namely
$$
     S^d_\mu(t) = \sum_{n= -\infty}^\infty   \chi_{I_n}(t) \frac{\mu(T_n)}{2^n}  =  
      \sum_{k=0}^\infty       \sum_{n=-\infty}^\infty  \chi_{I_n \ohne I_{n-1}}(t)     \frac{\mu(T_{n+k})}{2^{n+k}} =
      \sum_{k=0}^\infty S^d_{\mu, k}(t),
$$
where
$$
    S^d_{\mu, k}(t) =  \sum_{n=-\infty}^\infty  \chi_{I_n \ohne I_{n-1}}(t)     \frac{\mu(T_{n+k})}{2^{n+k}} = S^d_{\mu,0}(2^k t).
$$
Let us now look at an estimate from above for $S_\mu$ in terms of the $S^d_{\mu, k}$. Let $t \in I_n\ohne I_{n-1}$. Then
$$
   \max_{z \in T_k}  p_z(t) \le  \left\{ \begin{matrix}    \frac{2^k}{\pi 2^{2n}} &\text{ if } & n >  k ,\\
                                                                                              \frac{1}{\pi 2^{k}} &\text{ if } & n \le k. \\
                                                                                               \end{matrix} \right. 
$$
This easily implies that for $t \in I_n\ohne I_{n-1}$,
\begin{eqnarray*}
   S_{\mu}(t)   &\le&  \frac{1}{\pi} \left( \sum_{j=-\infty}^{-1}    \frac{2^{n+j}}{\pi 2^{2n}} \mu(T_{n+j})   +   \sum_{ j=0}^\infty         \frac{1}{2^{n+j}}   \mu(T_{n+j})   \right)    \\
                         &=&  \frac{1}{\pi} \left( \sum_{j=-\infty}^{-1}    2^{2j} \frac{\mu(T_{n+j})}{2^{n+j}}   +   \sum_{ j=0}^\infty        \frac{\mu(T_{n+j})}{2^{n+j}}   \right)    \\
                          &=&  \frac{1}{\pi} \left( \sum_{j=-\infty}^{-1}    2^{2j} S^d_{\mu,0}(2^j t)  +   \sum_{ j=0}^\infty         S^d_{\mu,0}(2^j t) \right).    \\
\end{eqnarray*}
Hence, for $t \in \RR$
\begin{eqnarray}    \label{eq:lowerbalay}
   S_\mu(t) &\le& \frac{1}{\pi} \left( \sum_{j=-\infty}^{-1}    2^{2j} S^d_{\mu,0}(2^j t)  +   \sum_{ j=0}^\infty         S^d_{\mu,0}(2^j t) \right) \\\nonumber
       & =& \frac{1}{\pi} \left( \sum_{j=-\infty}^{-1}    2^{2j} S^d_{\mu,0}(2^j t)  +   S^d_{\mu}( t) \right) .
\end{eqnarray}


\begin{theorem} \label{thm:pgqemb}
 Let $\mu$ be a positive regular Borel measure supported in a
  sector $S(\theta) \subset \CC_+$, $0 < \theta < \frac{\pi}{2}$ and let $1 \le q  <
 p < \infty$. Then the following are equivalent:
 \begin{enumerate}
 \item    The embedding
$$
     \LL: L^p(\RR_+) \rightarrow L^q(\CC_+, \mu), \quad f \mapsto \LL f,
$$
is well-defined and bounded.
\item The sequence $(2^{-n q/p'} \mu(S_{n}))$ is in $\ell^{p/(p-q)}(\ZZ)$.

\item The sequence $(       2^{n/p} \| \LL k_{2^n}\|_{L^q_\mu} ) $   is in $\ell^{qp/(p-q)}(\ZZ)$.
\end{enumerate}
If $p'<q$, then the above is also equivalent to
\begin{enumerate}
\setcounter{enumi}{3}
\item   $t^{q(2-p)/p}S_{\mu} \in L^{p/(p-q)}(\RR)$.
\end{enumerate}
\end{theorem}

{\bf Remark.} If $p'\ge q$, then the sweep $S_\mu$ may be infinite everywhere, so we cannot expect a characterisation in terms of $S_\mu$.

\proof
We will start by showing (2) $\Rightarrow$ (1).
Recall that $S(\theta)$ is contained in
$\bigcup_{n,k \in \ZZ,  |k| \le N }^\infty T_{n,k}$ for some $N \in \NN$. Suppose that (2) holds.
Now as in the proof of Theorem \ref{thm:pqemb}, we obtain for $z \in S_n$
$$
 |\LL f(z)| \le \int_0^\infty |e^{-zt}| |f(t)| dt \le \int_0^\infty
 |e^{-2^{n-1}t}| |f(t)| dt \le C_\Theta 2^{-n+1} Mf(2^{-n+1}).
$$
Note that $\chi_{I_{-n}}  Mf(2^{-n+1}) \le \chi_{I_{-n}} Mf$.
Consequently
\begin{eqnarray*}
  &&   \int_{\CC_+} |\LL f(z)|^q d\mu(z)  \\
&  \lesssim &    \sum_{n = -
     \infty}^\infty 2^{q(-n+1)} Mf(2^{-n+1})^q \mu(S_{n})\\
&   &    \sum_{n  = -
     \infty}^\infty 2^{q(-n+1)} Mf(2^{-n+1})^q \mu(S_n)\\
     &  \le &    \left( \sum_{n = -
     \infty}^\infty 2^{(-n+1)} Mf(2^{-n+1})^p  \right)^{q/p}   
     \left(    \sum_{n = -
     \infty}^\infty        2^{(q/p')(-n+1) (p/q)'}     \mu(S_{n})^{(p/q)'}  \right)^{1/(p/q)'}     \\    
& \lesssim &    \|f\|^{q}_{L^p}  \| (2^{-n q/p'} \mu(S_{n})) \|_{(p/q)'} \\
  \end{eqnarray*}
by the boundedness of the Hardy--Littlewood maximal function on $L^p(0, \infty)$.

(1) $\Rightarrow$ (2)  
For $\lambda>0$ we write $\tilde k_\lambda$ for the $L^p(0,\infty)$ function given by
\[
\tilde k_\lambda(t)=\lambda^{1/p}e^{-\lambda t} \qquad (t \ge 0),
\]
noting that $\|\tilde k_\lambda\|_{L^p}\asymp 1$.


From a result of Gurari\u\i\ and   Macaev in \cite{GM66}, it can be deduced that
\begin{equation}    \label{eq:macaev}
      \left\| \sum_{n \in \ZZ}  \alpha_n  \tilde k_{2^n} \right\|_p \approx \left( \sum_{n \in \ZZ} |\alpha_n|^p \right)^{1/p}
\end{equation}
for any $l^p$ sequence $(\alpha_n)$.

More precisely, we claim that there are constants $A,B>0$ such that for all scalars $(\alpha_k)$ we have
\[
A \sum_n |\alpha_n|^p \le \left\| \sum_n \alpha_n \tilde k_{2^n} \right\|_p^p \le B \sum_n |\alpha_n|^p
\]
By the change of variable $x=e^{-t}$ we have
\[
\left\| \sum_n \alpha_n \tilde k_{2^n} \right\|_p^p = \int_0^1 \left|\sum_n \alpha_n 2^{n/p} x^{2^n-1/p} 
\right|^p \, dx.
\]

Recall that a sequence $(m_j)$ is a {\em lacunary sequence\/} if $\inf m_{j+1}/m_j = r > 1$. 
Now, the result of Gurari\u\i\ and   Macaev in \cite{GM66} asserts the following:    If
$(n_j+1/p)$ is lacunary, then the sequence of functions $(t \mapsto (n_j+1/p)^{1/p}t^{n_j})$ in $L^p(0,1)$ is equivalent to the standard 
basis of $\ell^p$.

Writing $n_j=2^j-1/p$ for $j \in \ZZ$, we have the conditions of the  Gurari\u\i--Macaev theorem, and the claim follows.\\

%
%

Now,  denoting by $(\Omega, d\Omega)$ the probability space of sequences $(\epsilon_n)$ taking values in $-1$, $1$ with equal probability, equipped with the standard
product $\sigma$-algebra and probability measure, 
 we obtain
\begin{eqnarray*}
      \left( \sum_{n} |\alpha_n|^p \right)^{q/p} 
     & \gtrsim& \int_{\Omega}  \|  \LL( \sum_n \eps_n \alpha_n  \tilde k_n)\|^q_{L^q_\mu} d \Omega(\eps) \\
    &  =&   \int_{\CC_+}  \int_{\Omega} \left| \sum_n  \eps_n \alpha_n  2^{n/p} \frac{1}{2^n +z} \right|^q d\mu(z)  d \Omega(\eps) \\
    &\approx&    \int_{\CC_+}     \left(  \sum_n | \alpha_n|^2  2^{2n/p} \frac{1}{|2^n +z|^2} \right)^{q/2}  d\mu(z)  \\
     &=&   \sum_{k}     \int_{S_k}   \left(     \sum_n | \alpha_n|^2  2^{2n/p} \frac{1}{|2^n +z|^q} \right)^{q/2} d\mu(z)  \\
    &\gtrsim&      \sum_n \int_{S_n}  | \alpha_n|^q  2^{nq/p} \frac{1}{2^{qn}}    d\mu(z) \\
    & =&    \sum_n  | \alpha_n|^q 2^{-nq/p'} \mu(S_n).
      \end{eqnarray*}
Here, we have used the fact that $\mu$ is supported in a sector $S(\theta) $ in the last inequality.
Thus $(2^{-nq/p'} \mu(S_n))$ is a $l^{(p/q)'}$ sequence, and we have proved the desired implication.

A simple argument, again using sectoriality, shows that
$$
    \| \LL k_{2^n} \|_q^q \gtrsim \frac{1}{2^{nq}} \mu(S_n).
$$
Hence (3) $\Rightarrow$ (2). 

For  (1) $\Rightarrow$ (3), note that for any $n$,
$$
    \Re \LL  k_{2^n} \gtrsim | \LL  k_{2^n} (z)| \quad   \text{ for } z \in S(\theta).
$$
Hence for any sequence $(\alpha_n) \in l^{p/q}(\ZZ)$, $\alpha_n \ge 0$ for all $n$,
\begin{eqnarray*}
\sum_n \alpha_n   \| 2^{n/p}  \LL  k_{2^n} \|_{L^q(\mu)}^q
&\le & \int_{\CC_+}    \left(    \sum_n \alpha_n^{1/q}   2^{n/p}  |\LL  k_{2^n}(z) |   \right)^q d\mu(z)  \\
&\lesssim & \int_{\CC_+}    \left( \Re    \sum_n \alpha_n^{1/q}   2^{n/p}  \LL  k_{2^n}(z)    \right)^q d\mu(z)  \\
&\le & \int_{\CC_+}    \left| \LL  (\sum_n \alpha_n^{1/q}   2^{n/p}   k_{2^n}) (z)    \right|^q d\mu(z)  \\
&\lesssim & \left\|\sum_n \alpha_n^{1/q}   2^{n/p}   k_{2^n}  \right\|^{q}_p \\
& \lesssim&    \| \alpha_n^{1/q} \|^q_p =  \| \alpha_n \|_{p/q}
\end{eqnarray*}
by (\ref{eq:macaev}). This proves the equivalence of the first three statements.


(2) $\Rightarrow$ (4) Again, we can assume without loss that  $0<\theta < \pi/2$, in which
case $S(\theta) \subseteq \cup_{n = - \infty}^\infty T_n$.  If $t^{q(2-p)/p}S_{\mu} \in L^{p/p-q}(\RR)$, then by (\ref{eq:upperbalay}) $t^{q(2-p)/p}S^d_{\mu,0} \in L^{p/(p-q)}(\RR)$,
and
\begin{eqnarray*}
 &&    \sum_{n= -\infty}^\infty  2^{-nq p/(p'(p-q))} \mu(S_n)^{p/(p-q)}     \\   
    & = &   \sum_{n= -\infty}^\infty  2^n    2^{-n(p-2)q/(p-q)} \frac{\mu(S_n)^{p/(p-q)}}{2^{np/(p-q)}}      \\
    &    \approx  &
     \sum_{n= -\infty}^\infty   \int_{I_{n+1} \ohne I_n}  \left|    t^{-(p-2)q/p}S^d_{\mu,0} \right|^{p/(p-q)} dt \\
     &     = &      
     \|   t^{-(p-2)q/p}S^d_{\mu,0} \|^{p/(p-q)}_{p/(p-q)}                  < \infty.
\end{eqnarray*}
Thus (2) holds. 

Conversely, if $(2^{-n q/p'} \mu(S_{n})) \in l^{p/(p-q)}$, then   $t^{-(p-2)q/p}S^d_{\mu,0}    \in L^{p/p-q}(\RR)$ by the above calculation. By (\ref{eq:lowerbalay}),

\begin{eqnarray*}
 && \|   t^{-(p-2)q/p}S_{\mu} \|_{p/(p-q)}   \\
  &  \le& \frac{1}{\pi}   \left( \|   t^{-(p-2)q/p} S^d_{\mu}(t) \|_{p/(p-q)}   + \sum_{k= -\infty}^{-1} 2^{2k}  \|   t^{-(p-2)q/p}S_{\mu,0}^d(2^k t) \|_{p/(p-q)}   \right).\\
 \end{eqnarray*}

One sees easily that  
\begin{multline*}
 2^{2k}  \|   t^{-(p-2)q/p}S_{\mu,0}^d(2^k t) \|_{p/(p-q)}    =    2^{2k} 2^{-k(p-q)/p} 2^{k(p-2)q/p}  \|   t^{-(p-2)q/p} S_{\mu,0}^d \|_{p/(p-q)} \\
   = 2^{k (q +1 - q/p)} \|   t^{-(p-2)q/p} S_{\mu,0}^d \|_{p/(p-q)} ,
\end{multline*}  
thus the second term in the sum converges and is controlled by  the expression $ \|   t^{-(p-2)q/(p-q)} S_{\mu,0}^d \|_{p/(p-q)}$.
 For the first term, write
%
$$
  S^d_\mu(t) = \sum_n  \chi_{I_n}(t) \frac{\mu(T_n)}{2^n}  =   \sum_{k=0}^\infty       \sum_n \chi_{I_n \ohne I_{n-1}}(t)     \frac{\mu(T_{n+k})}{2^{n+k}}.
$$   
as before.   For each $k \ge 0$, it follows that
\begin{eqnarray*}
 &&  \int    \left(     t^{q(2-p)/p}    \sum_n   \chi_{I_n \ohne I_{n-1}}(t)     \frac{\mu(T_{n+k})}{2^{n+k}}                 \right)^{p/(p-q)} dt   \\
 &  \lesssim &     \sum_n 2^n  2^{nq(2-p)/(p-q)}  2^{-(n+k)p/(p-q)}\mu(T_{n+k})^{p/(p-q)} \\
 &  = &     \sum_n     2^{n (q/p -q  + 1)p/(p-q)}  2^{-(n+k)p/(p-q)} \mu(T_{n+k})^{p/(p-q)} \\
 &  = &     \sum_n     2^{n (-q/p'  + 1)p/(p-q)}  2^{-(n+k)p/(p-q)} \mu(T_{n+k})^{p/(p-q)} \\
 &  = &    2^{-k(1 -q/p'  )p/(p-q)} \sum_n     2^{(n+k) (-q/p'  )p/(p-q)}   \mu(T_{n+k})^{p/(p-q)} \\
 & \lesssim & 2^{-k(1 -q/p'  )p/(p-q)}.
\end{eqnarray*}
Hence $  t^{-(p-2)q/p}S_{\mu} \in L^{p/(p-q)}(\RR)$.   This concludes the proof.\qed

\subsection{A counterexample}
Let $\mu$ denote the measure on the interval $[1,\infty)$ defined by
$d\mu(x)=dx/\sqrt{x}$. Clearly, $\mu$ is sectorial and contained in a shifted
half-plane. Moreover $\mu$ satisfies the estimate that for a Carleson square $Q$
of size $h$ one has $\mu(Q) \le 2 h^{1/2}$.

Nonetheless, the Laplace--Carleson embedding $\LL: L^2(0,\infty) \to L^1(\mu)$ is unbounded (or equivalently,
in this case, the Carleson embedding $H^2(\CC_+) \to L^1(\mu)$ is unbounded). This can be
seen by noting that $\mu$ does not satisfy the condition of \cite[Theorem C]{luecking91}, since
the function $t \mapsto \int_{\Gamma(t)} x^{-1} d\mu(x)$ behaves as $x^{-1/2}$ and thus does not lie
in $L^2$. (Here $\Gamma(t)$ may be taken to be the interval $[t,\infty)$.)

It is constructive to give an explicit counterexample, following the reasoning of the proof of
\cite[Theorem C]{luecking91}. (Note that counterexamples in the case of the disc are simpler,
and can be found in \cite{TW72}.)

Define $\phi: \RR \to \RR$ by
\[
\phi(t)=\begin{cases}
1 & \text{if  $|t| \le 1$,}\\
t^{-1/2}(1+\log|t|)^{-1}  &\text{if $|t| \ge 1$.}
\end{cases}
\]
Thus $\phi \in L^2(\RR)$ and there is a function $F \in H^2$ with $\Re F(it)=\phi(t)$.
Now, if $F \in L^1(\mu)$ we would have
\[
\int \Re F \, d\mu \le  \int |F|  \, d\mu  < \infty,
\]
from which, as in \cite{luecking91} we could conclude (by writing $\Re F(x,0)$ in terms of the Poisson kernel) that
\[
A:=\int_{-\infty}^\infty \phi(t) \int_1^\infty \frac{x}{x^2+t^2 }\, \frac{dx}{\sqrt{x}} \, dt < \infty,
\]
and hence
\[
\int_{-\infty}^\infty  \phi(t)|t|^{-1/2} \, dt = 
\int_{-\infty}^\infty \phi(t) \int_{|t|}^\infty \frac{x}{2x^2 } \,\frac{dx}{\sqrt{x}} \, dt \le A < \infty,
\]
which is a contradiction.

\subsection{Measures supported in a strip}
\begin{theorem}
 \label{thm:pqexp1} Let $\mu$ be a positive regular Borel measure supported in a strip
 $\CC_{\alpha_1, \alpha_2} = \{ z \in \CC: \alpha_2 \ge \re z \ge  \alpha_1 \}$ for some $\alpha_2 \ge \alpha_1  >0$, and let $1 < p' \le q <\infty$, $q \ge 2$.
 Then the following are equivalent:
\begin{enumerate}
\item The embedding
$$
     \LL: L^p(0,\infty) \rightarrow L^q(\CC_+, \mu), \quad f \mapsto \LL f,
$$
is well-defined and bounded, with a bound only depending on the Carleson--Duren constant $C$ and the ratio $\frac{\alpha_2}{\alpha_1}$.
\item There exists a constant $C>0$ such that 
\begin{equation}  
\mu(Q_I) \le C
  |I|^{q/p'} \text{ for all intervals }   I \subset i\RR.
  \end{equation}
   \item  There exists a constant $C>0$ such that $\|\LL e^{- \cdot z}
  \|_{L^q_\mu} \le C \|e^{- \cdot z}\|_{L^p}$ for all $z \in \CC_{\alpha_1, \alpha_2}$.
  \end{enumerate}
\end{theorem}
\proof 
Again, obviously (1) $\Rightarrow$ (3). To show (3) $\Rightarrow$ (2), we  have to remember that the argument in Proposition \ref{prop} only works for Carleson squares $Q_I$ with centre 
$\lambda_I \in \CC_{\alpha_1, \alpha_2}$. Any Carleson square with centre 
$\lambda_I \in \CC_{\alpha_1/3, \alpha_2}$ can be covered by a Carleson square of at most triple sidelength with centre in
$\CC_{\alpha_1, \alpha_2}$, any Carleson square with centre in $\lambda_I \in \CC_{0,\alpha_1/3}$ has nonempty intersection with the support of $\mu$. 
If $Q_I$ is a Carleson square with centre  $\lambda_I \in \CC_{\alpha_2}$, then its intersection with the support of $\mu$ can be covered by at most
$\left[\frac{2 \re \lambda_I}{\alpha_2}   \right]+1$ Carleson squares with centre in $\CC_{\alpha_1, \alpha_2}$.
Hence
$$
   \mu(Q_I)   \lesssim \frac{|I|}{\alpha_2} \alpha_2^{q/p'}  \le |I|^{p/q'}.
$$

 This leaves (2) $\Rightarrow$ (1).

Consider the line parallel
to the imaginary axis $i\RR + \alpha_1/2$. Note that
$$ \LL: L^p(\RR_+)  \rightarrow L^2(i \RR + \frac{\alpha_1}{2})
$$
is bounded, since 
\begin{multline*}
   \|\LL f\|_{ L^2(i \RR + \frac{\alpha_1}{2})} = \|\LL (e^{-\alpha_1/2 t} f) \|_{L^2(i \RR )} = \|e^{-\alpha/2 t} f \|_2 \\
       \le  \|  e^{-\alpha_1 t}\|^{1/2}_{p/(p-2)}      \|f\|_p  \lesssim    \alpha_1^\frac{2-p}{2p}  \|f\|_p   .
\end{multline*}

Since the measure $\mu$ is supported in $\CC_{\alpha_1, \alpha_2} $, by the Carleson condition (\ref{eq:duren}) we have for each Carleson square in 
$Q_I = \{ z \in \CC:   i\Im z \in I, \alpha_1/2 < \Re z <  \alpha_1/2 +|I| \}$ in
$\CC_{+,\alpha_1/2} $:
$$
   \mu (Q_I) \le C |I|^{q/p'}    \le C    |I|^{q(1/p' -1/2) }   |I|^{q/2}   \le C    \alpha_2^{q(1/p' -1/2) }   |I|^{q/2}.
$$
Thus the Poisson embedding
$$
   L^2(i \RR + \frac{\alpha_1}{2}) \rightarrow  L^q(\CC_+, \mu)
$$
is bounded by Duren's Theorem, with constant $C\alpha_2^{1/p' -1/2 } $. Again, composing both maps yields the Laplace transform
$$
  \LL: L^p(\RR_+) \rightarrow L^q(\CC_+, \mu_\alpha)
$$
with norm bound $C\alpha_1^{\frac{2-p}{2p}} \alpha_2^{1/2 -1/p }  =  C \left(\frac{\alpha_2}{\alpha_1}\right)^{1/2 - 1/p}$.
\qed

\subsection{Sobolev spaces}

In this subsection, we will be interested in embeddings
$$
   \HH^2_\beta(0,\infty) \hookrightarrow L^p(\CC_+, \mu),
$$
where for $\beta>0$ the space  $\HH^2_\beta(0,\infty)$ is given by
\begin{eqnarray*}
 \HH^p_{\beta}(0,\infty) &=& \left\{ f \in L^p(\RR_+): \int_0^\infty
    |(\frac{d}{dx})^\beta f(t)|^p dt < \infty\right\},\\
    \|f\|_{\HH^p_\beta}^p &=&\|f\|^p_p + \| (\frac{d}{dx})^\beta f\|^p_p.
\end{eqnarray*}
Here $(\frac{d}{dx})^\beta f$ is defined as a fractional derivative
via the Fourier transform.

It is now easy to find versions of Theorems \ref{thm:pqemb} and
\ref{thm:pgqemb} for Sobolev spaces.
\begin{corollary} \label{thm:pqsobemb} Let $\mu$ be a positive Borel measure supported in a
  sector $S(\theta) \subset \CC_+$, $0 < \theta < \frac{\pi}{2}$, and
  let $q \ge p >1$. Then the following are equivalent:
\begin{enumerate}
\item The embedding
$$
     \LL: \HH_\beta^p(0,\infty) \rightarrow L^q(\CC_+, \mu), \quad f \mapsto \LL f,
$$
is well-defined and bounded.

\item There exists a constant $C>0$ such that $\mu_{q,\beta}(Q_I) \le C
  |I|^{q/p'}$ for all intervals in $I \subset i\RR$ which are symmetric
  about $0$. Here, $d\mu_{q,\beta}(z) = (1+ \frac{1}{|z|^{q \beta}}) d \mu(z)$.

\item  There exists a constant $C>0$ such that $\|\LL e^{- \cdot z}
  \|_{L^q_\mu} \le C \|e^{- \cdot z}\|_{\HH_\beta^p}$ for all $z \in \RR_+$.

\end{enumerate}
\end{corollary}
\proof
Follows immediately from Theorem \ref{thm:pqemb} and basic properties
of the Laplace transform.
\qed

\begin{corollary} \label{thm:pgqsobemb}
 Let $\mu$ be a positive regular Borel measure supported in a
  sector $S(\theta) \subset \CC_+$, $0 < \theta < \frac{\pi}{2}$ and let $1 \le q  <
 p$, $\beta \ge 0$. Suppose that $S_{\tilde \mu_{\beta,q}} \in L^{p/(p-q)}$. 
 Then the embedding
$$
     \LL: \HH_{\beta}^p(0,\infty) \rightarrow L^q(\CC_+, \mu), \quad f \mapsto \LL f,
$$
is well-defined and bounded.
\end{corollary}
\proof
Follows immediately from Theorem \ref{thm:pgqemb}.
\qed

 Laplace--Carleson embeddings of Sobolev spaces $\HH^2_\beta$ are
easily understood by means of
Theorem \ref{thm:ccarleson}:
 
\begin{theorem}
 \label{thm:sobolevcarl}
Let $\mu$ be a positive Borel measure on the  right half plane
$\CC_+$ and let $\beta >0$. Then the following are equivalent:
\begin{enumerate}
\item The Laplace--Carleson embedding
$$
     \HH^2_\beta(0,\infty) \rightarrow L^2( \CC_+, \mu)
$$
is bounded.
\item The measure 
$|1+ z|^{-2 \beta} d \mu(z)$ is a Carleson measure on $\CC_+$.
\end{enumerate}
\end{theorem}

\proof The proof is a simple reduction to the Carleson
  embedding theorem. Note that the map 
$$
    \HH^2_\beta(\RR_+) \rightarrow H^2(\CC_+) , \quad  f \mapsto (1+z)^{\beta} \LL f,
$$
is an isomorphism. The remainder follows from the holomorphy of $(1+z)^{\beta}$ and  $\LL f$ on $\CC_+$, and a density
argument.
\qed

\section{An application}
\label{sec:applic}

Suppose that we are given a
$C_0$--semigroup $(T(t))_{t \ge 0}$
defined on a Hilbert space $H$, with infinitesimal generator $A$, and
consider the system
\begin{equation}\label{eq:ABsystem}
{dx(t) \over dt}=Ax(t)+Bu(t), \qquad x(0)=x_0, \quad t \ge 0.
\end{equation}
Here $u(t) \in \mathbb C$ is the {\em input\/} at time $t$,
and 
$B: \mathbb C \to D(A^*)'$, the {\em control operator}.
We write $D(A^*)'$ for the completion of $H$ with respect to
the norm
\[
\|x\|_{D(A^*)'}=\|(\beta-A)^{-1}x\|_H,
\]
for any $\beta\in\rho(A)$.
To guarantee that the {\em state\/} $x(t)$ lies in $H$
one asks that $B \in \mathcal{L}(\mathbb C,D(A^*)')$ and  
\[
 \left\|\int_0^\infty T(t) Bu(t) \, dt \right\|_H \le m_0 \|u\|_{L^2(0,\infty)}
\]
for some $m_0>0$
(the admissibility condition for $B$). The $C_0$--semigroup $(T(t))_{t \ge 0}$ has an extension to $ D(A^*)'$.

We refer to the
survey \cite{jp} and the book \cite{TW} for the basic background to admissibility in the context of well-posed systems.
For diagonal semigroups, it is linked with the theory of Carleson measures as in
\cite{HR83,weiss88}: if $A$ has a Riesz basis of eigenvectors, with
eigenvalues $(\lambda_k)$, then a scalar control operator corresponding to a sequence $(b_k)$
is admissible if and only if the measure
\[
\mu:=\sum_k |b_k|^2 \delta_{-\lambda_k}
\]
is a Carleson measure for  $H^2(\CC_+)$.
An extension to normal semigroups
has also been made \cite{weiss99}.

Generalizations to $\alpha$-admissibility, in which
$u$ lies in $L^2(0,\infty; t^{\alpha} \, dt)$ for $-1<\alpha<0$, were studied by Wynn \cite{wynn}.
He used the fact that the Laplace transform maps
$L^2(0,\infty; t^{\alpha} \, dt)$ to a weighted Bergman space, for which a Carleson measure theorem is known.
The results above enable us to take this generalization further, considering admissibility
in the sense of the input lying in much more general spaces $L^2(0,\infty; w(t) \, dt)$.\\

Assume
that $1 \le q < \infty$ and the semigroup $(T(t))_{t \ge 0}$ acts on a Banach space $X$
with a $q$-Riesz basis of eigenvectors $(\phi_k)$; that is, $T(t)\phi_k=e^{\lambda_k t}\phi_k$ 
for each $k$, and
$(\phi_k)$ is a
a Schauder basis of $X$ such that
for some $C_1, C_2 > 0$ we have
\[
C_1 \sum |a_k|^q \le \| \sum a_k \phi_k \|^q \le C_2 \sum |a_k|^q
\]
for all sequences $(a_k)$ in $\ell^q$. 
Without loss of generality   $X=\ell^q$ and  
the eigenvectors of the generator of $(T(t))_{t \ge 0}$, denoted by $A$, are the canonical basis of $\ell^q$.
Suppose also that we have a Banach
space $Z$ of functions on $(0,\infty)$, either an $L^p$ space or a weighted space
$L^2(0,\infty; w(t)dt)$, whose dual space $Z^*$ can be regarded, respectively, as 
either $L^{p'}(0,\infty)$
or  $L^2(0,\infty; w(t)^{-1}dt)$ in a natural way. 

The following general result appears in \cite{jpp12}.

\begin{theorem} 
 Let $B$ be a linear bounded map from $\mathbb C$ to $D(A^*)'$ corresponding to the sequence $(b_k)$. The control operator $B$ is $Z$-admissible for $(T(t))_{t \ge 0}$, that is,
there is a constant $m_0>0$ such that
\[
 \left\|\int_0^\infty T(t) Bu(t) \, dt \right\|_X \le m_0 \|u\|_Z,\quad u\in Z,
\]
if and only if 
the Laplace transform induces a continuous mapping from $Z$
into $L^{q}(\CC_+,d\mu)$, where $\mu$ is the measure $\sum |b_k|^{q}\delta_{-\lambda_k}$.
\end{theorem}

This gives a direct application of our results on Laplace--Carleson embeddings: full details are given in \cite{jpp12}.


\section*{Acknowledgements}
This work was supported by EPSRC grant EP/I01621X/1.
The third author also acknowledges partial support of this work by a Heisenberg Fellowship of the German Research Foundation (DFG).


\begin{thebibliography}{99}










%

%

\bibitem{aleman} A. Aleman and A. Siskakis, Integration operators on Bergman spaces. {\em Indiana Univ. Math. J.} 46 (1997), no. 2, 337--356.

\bibitem{DP94}
N. Das and J.R. Partington,  
Little Hankel operators on the half-plane. 
{\em Integral Equations Operator Theory\/} 20 (1994), no. 3, 306--324. 


\bibitem{duren1}  P.L. Duren, Extension of a theorem of Carleson. \emph{Bull. Amer. Math. Soc.} 75 (1969) 143--146.

\bibitem{duren} P.~Duren, {\em Theory of $H^p$ spaces}. Academic Press, New York and London, 1970.

\bibitem{DGM}
P. Duren, E.A. Gallardo-Guti\'errez and A. Montes-Rodr\'\i guez, 
A Paley--Wiener theorem for Bergman spaces with application to invariant subspaces. 
{\em Bull. Lond. Math. Soc.} 39 (2007), no. 3, 459--466. 

 \bibitem{pau1}  P. Galanopoulos and J. Pau, Hankel  operators on large weighted Bergman spaces.
 {\em Ann. Acad.  Scient. Fenn. Math.}  37 (2012), 635 -- 648.

\bibitem{garnett} J.B. Garnett, {\em Bounded Analytic Functions}. Graduate Texts in Mathematics 236, Springer 2007.

\bibitem{GM66}
V. Gurari\u\i\ and V.I. Macaev,   Lacunary power sequences in spaces $C$ and $L\sb{p}$.  {\em Izv. Akad. Nauk SSSR Ser. Mat.} 30  (1966), 3--14. 

\bibitem{haak} B.H. Haak,  On the Carleson measure criterion in linear systems theory. 
{\em Complex Anal. Oper. Theory\/} 4 (2010), no. 2, 281--299.
  
\bibitem{zen09}
Z. Harper,   Boundedness of convolution operators and input-output maps between weighted spaces. 
{\em Complex Anal. Oper. Theory\/} 3 (2009), no. 1, 113--146. 

\bibitem{zen10}
Z. Harper,  Laplace transform representations and Paley--Wiener theorems for functions on vertical strips. 
{\em Doc. Math.} 15 (2010), 235--254.



\bibitem{HR83}
L.F.~Ho and D.L.~Russell,
Admissible input elements for systems in Hilbert space and a Carleson measure criterion,
{\em SIAM J. Control Optim.}  21  (1983),  no. 4, 614--640.  Erratum,
{\em SIAM J. Control Optim.}  21  (1983),  no. 6, 985--986.

\bibitem{jp}
B.~Jacob and J.R.~Partington, 
Admissibility of control and observation operators for semigroups: a survey.  
{\em Current trends in operator theory and its applications}, 199--221, 
Oper. Theory Adv. Appl., 149, Birkh\"auser, Basel, 2004. 

  
\bibitem{jpp12}
B. Jacob, J.R. Partington and S. Pott,
    Applications of Laplace--Carleson embeddings to admissibility and controllability. Submitted, 2012.
{\tt http://arxiv.org/abs/1203.2666}.

\bibitem{janson} S. Janson, 
Hankel operators between weighted Bergman spaces.
{\em Ark. Mat.} 26 (1988), no. 2, 205--219. 

\bibitem{rochberg}
P. Lin and R. Rochberg,  Hankel operators on the weighted Bergman spaces with exponential type weights. {\em  Integral Equations Operator Theory} 21 (1995), no. 4, 460 -- 483. 


\bibitem{l2}
D.H. Luecking,
Representation and duality in weighted spaces of analytic functions.
{\em Indiana Univ.~Math.~J.} 34 (1985), no. 2, 319--336.

\bibitem{luecking91}
D.H. Luecking, Embedding derivatives of Hardy spaces into Lebesgue spaces.
{\em Proc. London Math. Soc.} (3) 63 (1991), 595--619.



%

\bibitem{nik}
N.K.~Nikolski,   {\em Operators, functions, and systems: an easy reading. Vol. 1. Hardy, Hankel, and Toeplitz. 
Vol. 2. Model operators and systems.} 
Translated from the French by Andreas Hartmann and revised by the author. Mathematical Surveys and Monographs, 92--93. American Mathematical Society, Providence, RI, 2002.


\bibitem{ol}
V.L. Oleinik, Imbedding theorems for weighted classes of harmonic and analytic functions. 
 Investigations on linear operators and the theory of functions. \emph{V. Zap. Nauc. Sem. Leningrad. Otdel. Mat. Inst. Steklov.} 
(LOMI) 47 (1974), 120--137, 187, 192--193.
 Translated in: 
 \emph{J. Soviet Math.} 9(1978), no 1, 228--243.




%


\bibitem{pau}  J. Pau and J.A. Pel{\'a}ez, Embedding theorems and integration operators on Bergman spaces with rapidly decreasing weights. {\em J. Funct. Anal.}
 259 (2010), no. 10, 2727--2756.
 
 




 
 
 \bibitem{stein} E. Stein,
 Harmonic analysis: real-variable methods, orthogonality, and oscillatory integrals. 
 With the assistance of Timothy S. Murphy. Princeton Mathematical Series, 43. Monographs in Harmonic Analysis, III. Princeton University Press, Princeton, NJ, 1993.

\bibitem{ST80}
J.O. Str\"omberg and A. Torchinsky,  
Weights, sharp maximal functions and Hardy spaces.
{\em Bull. Amer. Math. Soc.} (N.S.) 3 (1980), no. 3, 1053--1056. 



%

%


%
%


%
%
%

\bibitem{TW72}
B.A. Taylor and D.L. Williams,   Interpolation of $l\sp{q}$ sequences by $H\sp{p}$ functions. 
{\em Proc. Amer. Math. Soc.} 34  (1972), 181--186.


\bibitem{TW}
M. Tucsnak and G. Weiss, {\em Observation and control for operator semigroups}. Birkh\"auser Advanced Texts: Basler Lehrb\"ucher.  Birkh\"auser Verlag, Basel, 2009. 

\bibitem{weiss88}
G. Weiss,  Admissibility of input elements for diagonal semigroups on $l\sp 2$,
{\em Systems Control Lett.} 10  (1988),  no. 1, 79--82.

\bibitem{weiss99}
G. Weiss, A powerful generalization of the Carleson measure theorem? {\em Open problems in mathematical systems and control theory}, 267--272, Comm. Control Engrg. Ser., Springer, London, 1999.

\bibitem{wynn}
A. Wynn,   $\alpha$-admissibility of observation operators in discrete and continuous time. 
{\em Complex Anal. Oper. Theory\/} 4 (2010), no. 1, 109--131. 

\bibitem{yong} Liu Yongmin, 
Small Hankel operators on weighted Bergman spaces of bounded symmetric domains. 
{\em Acta Math. Sci. Ser. B Engl. Ed.} 20 (2000), no. 1, 27--34.

\bibitem{zhu}  Ke He Zhu, 
 Operator theory in function spaces. Second edition. Mathematical Surveys and Monographs, 138. American Mathematical Society, Providence, RI, 2007. xvi+348 pp. 

\bibitem{zhu1} Ke He Zhu, 
Hankel operators on the Bergman space of bounded symmetric domains. 
{\em Trans. Amer. Math. Soc. } 324 (1991), no. 2, 707--730.

\end{thebibliography}
\end{document}